\documentclass[11pt]{article}      
\usepackage[dvips]{graphics}      
\usepackage{amssymb}      
\usepackage{amsmath}      
\usepackage{amsthm}      
\usepackage{amsbsy}      
\usepackage[dvips]{graphics}      
\usepackage{amscd}      
 
\theoremstyle{plain}      
\newtheorem{theorem}{Theorem}[section]

\newtheorem{lemma}[theorem]{Lemma}      
\newtheorem{corollary}[theorem]{Corollary}      
\newtheorem{proposition}[theorem]{Proposition}

\newtheorem{definition}{Definition}[section]      
\theoremstyle{remark}      
\newtheorem{remark}[theorem]{Remark}

\newcommand{\Z}{{\mathbb{Z}}}      
      
\newcommand{\R}{{\mathbb{R}}}

\usepackage{amssymb}      
\newcommand{\aps}{{{\alpha_P^{\star}}}}        
\newcommand{\bps}{{{\beta_P^{\star}}}}

\begin{document}

\title{The braided Ptolemy-Thompson group $T^*$ is asynchronously 
combable\footnote{First version: {\tt June 20, 2005}. This version:
 {\tt February 7, 2006}      
This preprint is available electronically at       
\tt  http://www-fourier.ujf-grenoble.fr/\~{ }funar } 
}      

\author{
\begin{tabular}{cc}      
 Louis Funar &  Christophe Kapoudjian\\      
\small \em Institut Fourier BP 74, UMR 5582       
&\small \em Laboratoire Emile Picard, UMR 5580\\      
\small \em University of Grenoble I &\small \em University of Toulouse      
III\\      
\small \em 38402 Saint-Martin-d'H\`eres cedex, France      
&\small \em 31062 Toulouse cedex 4, France\\      
\small \em e-mail: {\tt funar@fourier.ujf-grenoble.fr}      
& \small \em e-mail: {\tt ckapoudj@picard.ups-tlse.fr} \\      
\end{tabular}      
}

\maketitle      
  
\begin{abstract} 
The braided Ptolemy-Thompson group $T^*$ is an extension  
of the  Thompson group $T$ by the 
full braid group $B_{\infty}$ on infinitely many strands. This group is 
a simplified version of the acyclic extension considered 
by Greenberg and Sergiescu, and can be viewed as a   
mapping class group of a certain infinite planar surface. 
In a previous paper we showed that $T^*$ is finitely presented. 
Our main result  here is that $T^*$  (and $T$) is asynchronously combable. 
The method of proof is inspired by Lee Mosher's proof of automaticity 
of mapping class groups. 

\vspace{0.2cm} 
 
 
 
\end{abstract} 
 
{\small
\tableofcontents 
}     
\section{Introduction}   
\subsection{Statements and results}   
The Thompson groups $T$ and $V$ were the first examples 
of finitely presented infinite simple groups. We refer to \cite{CFP} 
for a survey concerning some of their properties. 
These  groups arise geometrically as groups of 
almost-automorphisms of the  infinite binary tree, where 
by an almost-automorphisms is meant an automorphism outside a bounded subset.

\vspace{0.2cm}
\noindent 
The group $V$ is a group of homeomorphisms of the Cantor set, and as such, 
it might be thought of as a group of infinite permutations. 
There is a well-known relation between  permutations and 
braids, in which one replaces transpositions  
by the usual braid generators. Moreover,  we can associate 
in the same way an Artin group to any Coxeter group. 
Unfortunately, this algebraic point of view seems to be too 
narrow in order to define a convenient analog of the Artin group 
for the group $V$. 

\vspace{0.2cm}
\noindent In searching for a universal mapping class group we found 
that a geometric  method (similar to the Artinification above) 
yields an interesting group, as follows. Consider first the surface obtained  
by thickening the binary tree and remark then  that 
permutations can be lifted to classes of  homeomorphisms  
of this surface which rather braid the boundary components 
instead to permute them. Lifting all elements of $V$ one finds 
the group ${\mathcal B}$, which we 
proved in  \cite{FK} that it is finitely presented.
This group can be viewed as the asymptotic mapping class group 
of a sphere minus a Cantor set. 
Soon afterwards,  M.Brin (\cite{Bri,Bri2}) and P.Dehornoy 
(\cite{De1,De2}) constructed  and studied groups that correspond to the 
asymptotic mapping class group of a disk minus a Cantor set, and which 
are therefore braid-like.  

\vspace{0.2cm}
\noindent
All these  groups are extensions of the group $V$ by an infinite 
{\em pure} mapping class (or braid) group.  
However, it has been shown in \cite{FK2} that we can do better 
if we restrict ourselves to the smaller group $T$. In fact, we constructed 
an extension $T^*$ of the group $T$ by the whole infinite 
braid group $B_{\infty}$. The group $T$ received a lot of attention 
since E.Ghys and V.Sergiescu  (\cite{GS}) 
proved that $T$ can be embedded in the 
diffeomorphism group of the circle and it can be viewed as a sort of 
discrete analog of the later. Moreover, the group $T^*$ is still an 
asymptotic mapping class group of a suitable surface of infinite type, 
which is homeomorphic to a thick tree whose edges are punctured. 
We found also that $T^*$ is a simplified version of the mysterious 
acyclic extension constructed by P.Greenberg and V.Sergiescu (\cite{GrS}). 
Our main result in \cite{FK2} is that $T^*$ is finitely presented and this 
was also a way to approach the finite presentability of the 
universal mapping class group of infinite genus. 

\vspace{0.2cm}
\noindent 
The aim of the present paper is to show that $T^*$ has strong finiteness 
properties. Although it was known  that we can generate the 
Thompson groups  using automata (\cite{GN}) very little was known about the 
geometry of their Cayley graph. Recently, V.Guba made progress on this 
question (\cite{Gu,GuS}). We wanted to approach this problem from the 
perspective of the mapping class groups, since we can view $T$ as a 
mapping class group of a surface of infinite type. One of the far reaching 
results in this respect is the Lee Mosher theorem  (\cite{Mo}) stating that 
mapping class of finite surfaces are automatic. Our main result 
shows that, when shifting to infinite surfaces, a slightly 
weaker result still holds true, namely: 

\begin{theorem}
The group $T^*$ is asynchronously combable. 
\end{theorem}
 
\noindent In particular, in the course of the proof we prove also that:

\begin{corollary}
The Thompson group $T$ is asynchronously combable. 
\end{corollary}

\noindent The proof was greatly inspired by the methods of L.Mosher. 
The mapping class group was embedded in the Ptolemy groupoid of some 
triangulation of the surface, as defined by L.Mosher and R.Penner.  
It suffices then to provide combings for the later. 

\vspace{0.2cm}
\noindent 
In our case the respective Ptolemy groupoid is, fortunately,  
the group $T$, which could be viewed also as a groupoid acting 
on triangulations of the hyperbolic plane. 
The first  difficulty consists of dealing with the fact that 
the surface under consideration is 
non-compact. Thus we have to get  extra control on the action of 
$T$ on triangulations and in particular to consider 
a finite set of generators 
of $T$ instead of  the set of all flips that was used 
by Mosher for compact surfaces.    
The second difficulty is that we should modify the Mosher algorithm 
in order to obtain the boundedness  of the combing. 
Finally, by  shifting from $T$  to $T^*$  it amounts to 
consider triangulations of the hyperbolic plane whose edges are punctured. 
The same procedure works also in this situation, but we need  
another ingredients to get explicit control on the braiding, which 
remind us the geometric  solution of the word problem for braid groups.

\vspace{0.2cm}
\noindent 
\noindent {\bf Acknowledgements.} The authors are indebted to  
Vlad Sergiescu and Bert Wiest for comments and useful discussions.

\subsection{Definition of  the braided Thompson group $T^*$}      
The main step in obtaining the universal mapping class group of genus zero  
${\mathcal B}$ is to replace compact surfaces by an infinite surface 
and to consider the mapping classes of those homeomorphisms 
having a nice behaviour at infinity.      
We have however to take into account various versions of the surface  
considered in \cite{FK}, notably the pointed and holed spherical 
or planar surfaces.

\vspace{0.2cm}
\noindent 
Let ${\cal T}$ be the infinite binary tree.  
A {\em finite binary tree} $T$ is a finite subtree of ${\cal T}$ whose      
  internal vertices are all 3-valent. Its terminal vertices (or 1-valent      
  vertices) are called {\em leaves}. We denote by ${\cal L}(T)$ the set of      
  leaves of $T$, and call the number of leaves the {\em level} of $T$.      
 
\begin{definition}[Thompson's group $V$]      
A symbol $(T_1,T_0,\sigma)$ is a triple consisting of two finite      
  binary trees $T_0$, $T_1$ of the same level, together with a bijection $\sigma:{ \cal L}(T_0)\rightarrow      
  { \cal L}(T_1)$.      
 
\vspace{0.2cm} 
\noindent 
If $T$ is a finite binary subtree of ${\cal T}$ and $v$ is a leaf of      
  $T$, we define the finite binary subtree $\partial_v T$ as the union of $T$      
  with the two edges which are the descendants of $v$. Viewing $\partial_v T$   
  as a subtree of the planar tree ${\cal T}$, we may distinguish the left descendant from the right descendant      
  of $v$. Accordingly, we denote by $v_l$ and $v_r$ the leaves of the two new      
  edges of $\partial_v T$.       
 
\vspace{0.2cm} 
\noindent 
  Let ${\cal R}$ be the equivalence relation on the set of symbols      
  generated by the following relations:      
$$(T_1,T_0,\sigma)\sim_v (\partial_{\sigma(v)}T_1,\partial_v T_0,\partial_v      
\sigma)$$      
where $v$ is any leaf of $T_0$, and $\partial_v \sigma: { \cal L}(T_0)\rightarrow      
  { \cal L}(T_1)$ is the natural      
extension of $\sigma$ to a bijection ${\cal L}(\partial_v T_0)\rightarrow {\cal L}(\partial_{\sigma(v)} T_1)$ which maps $v_l$ (resp. $v_r$) to $\sigma(v)_l$      
(resp. $\sigma(v)_r$). Denote by $[T_1,T_0,\sigma]$ the class of a      
symbol $(T_1,T_0,\sigma)$, and by $V$ the set of equivalence classes for      
the relation ${\cal R}$. Given two elements of $V$, we may represent them by      
two symbols of the form $(T_1,T_0,\sigma)$ and $(T_2,T_1,\tau)$ respectively, and define      
the product      
$$[T_2,T_1,\tau]\cdot [T_1,T_0,\sigma]=[T_2,T_0,\tau\circ\sigma]$$      
This endows $V$ with a group structure, with neutral element $[T,T,id_{{\cal      
    L}(T)}]$, where $T$ is any finite binary subtree. This is Thompson's group $V$ (cf. \cite{CFP}).         
\end{definition}      
     
\begin{definition}[Ptolemy-Thompson's group $T$]\label{tho}     
Let $T_{S_0}=\overline{\lambda(S_0)}$ be the smallest finite binary subtree     
  of ${\cal T}$ containing $\lambda(S_0)$. Choose a cyclic 
counterclockwise     
  labeling of its leaves by $1,2,3$. Extend inductively this cyclic labeling     
  to a cyclic labeling of the leaves of any finite binary subtree of ${\cal T}$     
  containing $T_{S_0}$: if $T_{S_0}\subset T\subset \partial_v T$, where $v$     
  is a leaf of a cyclically labeled tree $T$, then there is a unique cyclic     
  labeling of the leaves of $\partial_v T$ such that:     
\begin{itemize}     
\item if $v$ is not the leaf 1 of $T$, then it is also the leaf 1 of     
  $\partial_v T$;     
\item if $v$ is the leaf 1 of $T$, then the leaf 1 of $\partial_v T$ is the     
  left descendant of $v$.     
\end{itemize}     
\vspace{0.2cm} 
\noindent Thompson's group $T$ (also called Ptolemy-Thompson's group) is the subgroup of $V$ consisting of elements     
  represented by symbols $(T_1,T_0,\sigma)$, where $T_1,T_0$ contain     
  $\lambda(S_0)$, and $\sigma: { \cal L}(T_0)\rightarrow  { \cal L}(T_1)$ is a cyclic permutation. The cyclicity of   
  $\sigma$ means that there exists some integer $i_0$, $1\leq i_0\leq n$ (if   
  $n$ is the level of $T_0$ and $T_1$), such that $\sigma$ maps the $i^{th}$   
  leaf of $T_0$ onto the $(i+i_0)^{th}$ (mod $n$) leaf of $T_1$, for $i=1,\ldots, n$.     
\end{definition}

\vspace{0.2cm} 
\noindent
The surfaces below will be oriented and all      
homeomorphisms considered in the sequel will be      
orientation-preserving, unless the opposite is explicitly stated.   
\begin{definition}\label{ss}      
The ribbon tree $D$ is the planar surface obtained by   
thickening in the plane the infinite binary        
tree. We denote by $D^{\star}$ the ribbon tree  
with infinitely many punctures, one puncture for each edge of the 
tree. 
\end{definition}

\begin{figure}      
\hspace{5cm}\includegraphics{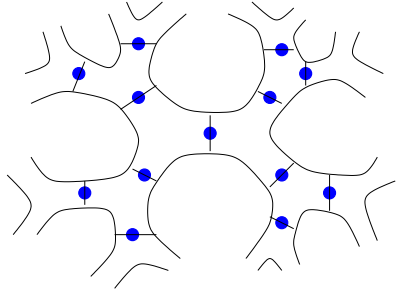}      
\caption{$D^{\star}$ and its canonical rigid structure}\label{tt}      
\end{figure}

\begin{definition} 
A {\em rigid structure} on $D^{\star}$ 
 is a decomposition  into punctured 
hexagons by means of a family of arcs through the punctures, whose 
endpoints are on the boundary of $D$. It is assumed that these arcs  
are pairwise non-homotopic in $D$, by homotopies keeping the boundary points 
on the boundary of $D$. There exists a canonical rigid structure,  
in which arcs are segments transversal to the edges, as drawn in the  
picture \ref{tt}.

\vspace{0.2cm} 
\noindent A planar subsurface of $D$ is {\em admissible}  if  
it is an union of hexagons coming from the canonical rigid structure.  
The {\em frontier} of an admissible surface is the  
union of the arcs contained in the boundary.    
\end{definition} 
 
\begin{remark} 
Two different arcs associated to the same puncture should be isotopic  
in $D$. Nevertheless, these arcs might be  
non-isotopic (non-homotopic) in $D^{\star}$ i.e.  if one asks that 
the isotopies keep fixed the punctures.  
\end{remark}

\begin{definition}      
Let $\varphi$ be a homeomorphism of $D^{\star}$. One says that       
$\varphi$ is {\em asymptotically rigid} if the following conditions are      
fulfilled:      
\begin{itemize}      
\item There exists an admissible subsurface $\Sigma\subset 
  D^{\star}$   such      
  that $\varphi(\Sigma)$ is also admissible.       
\item The complement $D^{\star} -\Sigma$  is an union of $n$ infinite      
  surfaces. Then the restriction  
$\varphi: D^{\star} -\Sigma\to D^{\star}-\varphi(\Sigma)$      
 is {\em rigid}, meaning that it respects the rigid structures in the      
 complementary of the compact subsurfaces, it  maps the 
hexagons into hexagons.       
  Such a surface $\Sigma$ is called a {\em support} for $\varphi$.      
\end{itemize}      
One denotes by $T^{\star}$ the 
{\em group of asymptotically rigid     homeomorphisms of $D^{\star}$} 
modulo isotopy through homeomorphisms  
which pointwise preserve the boundary $D^{\star}$. 
\end{definition} 
\begin{remark} 
There exists a cyclic order on the frontier arcs of an  
admissible subsurface induced by the planarity. An asymptotically 
rigid homeomorphism necessarily preserves the cyclic order  
of the frontier for any admissible subsurface.  
In particular one can identify $T$ with the group of asymptotically 
rigid homeomorphisms (mod isotopy) of the ribbon tree $D$.  
Further  $T^{\star}$  is the analogue of  
$T$  for the punctured  disk.  
\end{remark}

\subsection{Preliminaries on combings}
We will follow below the terminology introduced by Bridson in \cite{Br,Br2},
in particular we allow very general combings. We refer  the reader 
to \cite{ECHLPT} for a thorough introduction to the subject. 

\vspace{0.2cm}
\noindent 
Let $G$ be a finitely generated group with a finite generating set 
$S$, such that $S$ is closed with respect to the inverse 
and corresponding Cayley graph $C(G,S)$. This graph is endowed with 
the word metric in which the distance $d(g,g')$ between the vertices 
associated 
to the elements $g$ and $g'$ of $G$ is the minimal length of a word 
in the generators $S$  representing  the 
element $g^{-1}g'$ of $G$. 

\vspace{0.2cm}
\noindent 
A {\em combing} of the group $G$ with generating set $S$ 
is a map which associates to any element $g\in G$ a 
path $\sigma_{g}$ in the Cayley graph associated to $S$
from $1$ to $g$. In other words $\sigma_{g}$ is a word in the 
free group generated by $S$ that represents  the 
element $g$ in $G$.  We can also represent $\sigma_g(t)$ as a 
combing path in $C(G,S)$ that joins the identity element to $g$, moving 
at each step to a neighboring vertex and which becomes eventually 
stationary at $g$. Denote by $|\sigma_g|$ the length of the path $\sigma_g$ 
i.e. the smallest $t$ for which $\sigma_g(t)$ becomes stationary.  

\begin{definition}
The combing $\sigma$ of the group $G$ is {\em synchronously bounded} if it 
satisfies the synchronous fellow traveler property below. This means 
that there exists $K$ such that the combing paths $\sigma_g$ and 
$\sigma_{g'}$ of two elements $g$, $g'$ at distance $d(g,g')=1$ 
are at most distance $K$ far apart at each step i.e. 
\[ d(\sigma_g(t),\sigma_{g'}(t))\leq K, \; {\rm for \; any }\; t\in R_+\]
A group  $G$ having an asynchronously bounded combing is asynchronously
combable. 
\end{definition}

\vspace{0.2cm}
\noindent In particular, combings furnish normal forms for 
group elements. 
The existence of combings with special properties (like 
the fellow traveler property) has important consequences for the 
geometry of the group (see \cite{A,Br}).

\vspace{0.2cm}
\noindent 
We will introduce also a slightly weaker condition (after 
Bridson and Gersten) as follows:

\begin{definition}
The combing $\sigma$ of the group $G$ is {\em asynchronously bounded} if it 
satisfies the asynchronous fellow traveler property below. This means 
that there exists $K$ such that  for any 
two elements $g$, $g'$ at distance $d(g,g')=1$ there exist ways to 
travel through  the combing paths $\sigma_g$ and 
$\sigma_{g'}$ at possibly different speeds so that corresponding  
points are at most distance $K$ far apart.  Thus, there exists 
continuous increasing functions $\varphi(t)$ and $\varphi'(t)$ 
going from zero to infinity such that   
\[ d(\sigma_g(\varphi(t)),\sigma_{g'}(\varphi'(t)))\leq K, \; {\rm for \; any }\; t\in R_+\]
The asynchronously bounded combing $\sigma$ has {\em a departure function} 
$D:\R_+\to \R_+$ if, for all $r >0$, $g\in G$ and 
$0\leq s,t\leq |\sigma_g|$, the assumption  $|s-t| > D(r)$ implies that 
$d(\sigma_g(s),\sigma_g(t)) > r$. 
\end{definition}

\begin{remark}
There are known examples of asynchronously combable groups with 
a departure function: 
asynchronously automatic groups (see \cite{ECHLPT}), 
the fundamental group of a Haken 3-manifold (\cite{Br}), 
or of a geometric 3-manifold (\cite{Br2}),
semi-direct products of $\Z^n$ by $\Z$ (\cite{Br}). 
Gersten (\cite{Ger}) proved that such groups are of type ${\rm FP}_3$ and 
announced that they should actually be ${\rm FP}_{\infty}$.
Notice that there exist asynchronously combable groups  
(with departure function) which are 
not asynchronously automatic, for instance the {\sf Sol} and 
{\sf Nil} geometry groups  of closed 3-manifolds 
(see \cite{Bra});  in particular, they are not automatic. 
\end{remark}

\section{The Thompson group $T$ is asynchronously combable}
\subsection{The group $T$ as a mapping class group}
We consider the following elements of $T$, defined as mapping classes 
of asymptotically rigid homeomorphisms: 
\begin{itemize} 
\item The support of the element $\beta$ is the central hexagon.  
Further $\beta$ acts as the counterclockwise rotation of order  
three whose axis is vertical and which permutes the  three branches 
of the ribbon tree issued from the hexagon.  

\begin{center} 
\includegraphics{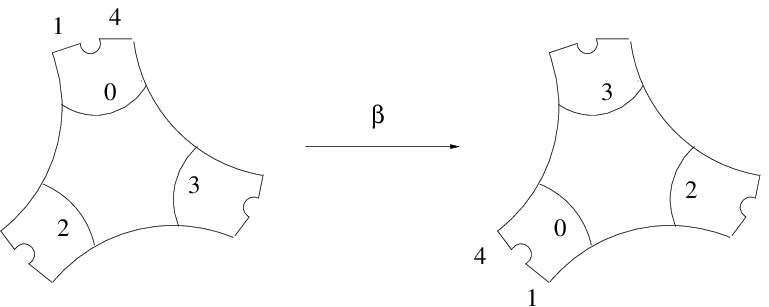}      
\end{center}
\item  The support of $\alpha$  is the union of two adjacent hexagons, 
one of them being the support of $\beta$.  Then $\alpha$  
rotates counterclockwise the support of angle $\frac{\pi}{2}$, by  
permuting the four branches of the ribbon tree
issued from the support.  

\begin{center} 
\includegraphics{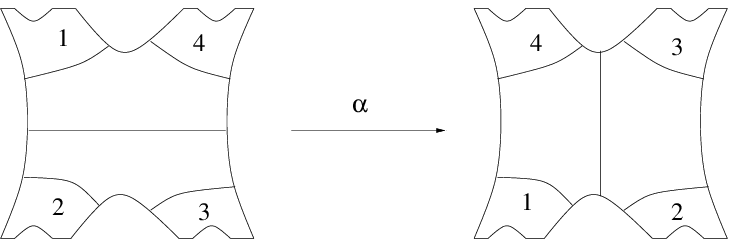}      
\end{center}
\end{itemize} 

\noindent Lochak and Schneps (\cite{LS}) proved 
 that the group $T$ has the following presentation with generators 
$\alpha$ and $\beta$ and relations 
\[ \alpha^4=\beta^3=1\]
\[ [\beta\alpha\beta, \alpha^2\beta\alpha\beta\alpha^2]=1\]
\[  [\beta\alpha\beta, \alpha^2\beta^2\alpha^2\beta\alpha\beta\alpha^2
\beta\alpha^2]=1\]
\[ (\beta\alpha)^5=1\]

\begin{remark}
If we set $A=\beta\alpha^2$, $B=\beta^2\alpha$ and 
$C=\beta^2$ then we obtain the generators $A,B,C$ of the group $T$, 
considered in \cite{CFP}. Then the  first two relations above 
are equivalent to 
\[ [AB^{-1}, A^{-1}B A]=1, \,\;\; [AB^{-1}, A^{-2}B A^2]=1\]
The presentation of $T$ in terms of the generators 
$A,B,C$ consists  of the two relations above with  four more 
relations to be added: 
\[ C^3=1, \;\, C=B A^{-1}CB,\,\; CA=(A^{-1}CB)^2,\,\; (A^{-1}CB)(A^{-1}BA)=B
(A^{-2}CB^2)\]
\end{remark}

\subsection{The Ptolemy groupoid and $T$}
The previous results  from \cite{FK2} come from the interpretation 
of the group $T$ (and its braided version $T^*$) 
as a mapping class group of an infinite surface. 
In this sequel we will bring forth another perspective,   
by turning back to Penner's  original approach (\cite{pe0,pe}) of 
the Ptolemy groupoid acting on triangulations of surfaces. 
In the case when the surface is the hyperbolic plane $H^2$ Penner 
obtained what he called the universal Ptolemy groupoid.
It was remarked  that this groupoid is actually a group since 
all elements have inverses.  
A result of Kontsevich and Imbert (\cite{im}) 
identified the Ptolemy groupoid with  both the group $T$ 
and $PPSL(2,\Z)$, as well. An explicit and 
convenient description of this isomorphism was described 
by Lochak and Schneps in \cite{LS}, and that permitted them  
to obtain a nice geometric-oriented  presentation of $T$. 

\vspace{0.2cm} 
\noindent  
Let us recall a few definitions which will be needed in the sequel. 
Details may be found in \cite{pe0,pe}. 
One considers ideal triangulations of the hyperbolic plane 
$H^2$ having vertices at the rational points of the boundary circle 
and coinciding with the Farey tessellation for all but finitely 
many triangles. Moreover, each triangulation is labeled, 
by choosing a distinguished oriented edge $e$ (called the 
d.o.e.). This theory has been considered for arbitrary surfaces 
in \cite{pe}. 

\vspace{0.2cm} 
\noindent  
Let $\gamma$ be an edge (i.e. an ideal arc) of the triangulation 
$\tau$ (unlabeled for the moment). Then $\gamma$ is a diagonal 
of a unique quadrilateral $Q\subset \tau$. Let $\gamma^*$ be 
the other diagonal of $Q$. The triangulation $\left(\tau -\{\gamma\}\right) \cup
\{\gamma^*\}$, obtained from $\tau$ by removing 
$\gamma$ and replacing it by $\gamma^*$ 
is said to be the result of applying the flip on the edge 
$\gamma$. We denote by $F_{\gamma}$ this flip. 
The definition extends to the labeled case without modifications 
when $\gamma$ is not the d.o.e., by keeping the same d.o.e.
When $\gamma$ is the d.o.e. we give the flipped triangulation the 
d.o.e. $\gamma^*$ with that orientation which makes the frame 
$\{\gamma,\gamma^*\}$ (in this order) positively oriented. 

\vspace{0.2cm} 
\noindent 
Flips can be composed and they form a (free) groupoid. 
Moreover, one assumes that two compositions are 
identical if their actions on a fixed  triangulation are identical. 
One obtains a  quotient groupoid, which is called 
the Ptolemy groupoid $PT$, and one verifies easily that 
all elements have inverses and thus this is actually a group. 
The group $PT$  is therefore generated by all  possible flips 
on edges of the triangulation.  Notice that there are 
infinitely many  different flips on the infinite triangulation. 
However, one can infer from \cite{im,LS} that the Ptolemy group 
$PT$ is just another instance of our familiar Thompson group $T$
(or the group of piecewise projective transformations $PPSL(2,\Z)$). 
In particular, $PT$ is finitely presented. 

\vspace{0.2cm} 
\noindent 
Let us explain now some details concerning the identification 
of the Ptolemy groupoid appearing in Lochak-Schneps picture with 
that considered by the present authors (see also in \cite{FK,KS}). 
Lochak and Schneps defined two generators of $PT$, which are  
the two  {\em local} moves below: 
\begin{itemize}
\item The  fundamental flip, which is the flip $F=F_e$ on the d.o.e. $e$. 
\item The rotation $R$ which preserves the triangulation but 
moves the given d.o.e. $e$ in  the clockwise direction to 
the next edge (adjacent to $e$) of the triangle sitting on the left 
of the d.o.e. and containing the d.o.e. as an edge.  
\end{itemize}
We wish to emphasize that these two moves are {\em local}. 
All other edges of the triangulation are kept pointwise fixed. 
It is not so difficult to show that the two local moves above generate 
the group $PT$, because an arbitrary flip can be obtained by 
conjugating the fundamental flip $F$ by a composition of 
rotations $R$ and orientation-reversals $F^2$ of the d.o.e. 

\vspace{0.2cm} 
\noindent 
There exists another way to look at the group $PT$, which makes the 
identification with $T$ manifest. An element of $PT$ is specified
by a couple of two labeled triangulations $(\Delta, \Delta')$ as above. 
We associate a homeomorphism of the closed disk 
$\overline{H^2}$ obtained by compactifying the open disk model 
of the hyperbolic plane,  which is subject to the 
following requirements: 
\begin{itemize}
\item The homeomorphism is piecewise linear with respect to the 
triangulations $\Delta$ and $\Delta'$. This means that it sends 
each triangle of $\Delta$  onto some triangle of $\Delta'$  
by a transformation from $PSL(2,\Z)$.  
\item The homeomorphism sends the d.o.e. of $\Delta$ onto the d.o.e. 
of $\Delta'$ with the corresponding orientation. 
\end{itemize}
The homeomorphism is then uniquely determined 
by the two conditions  above and it 
is an element of $PPSL(2,\Z)$.  It is also determined 
by its restriction to the boundary, when   
$PPSL(2,\Z)$ is viewed as a subgroup of ${\rm Homeo}_+(S^1)$. 
Denote by $\Phi:PPSL(2,\Z)\to PT$ the inverse correspondence. 
Recall that $PPSL(2,\Z)$ is isomorphic to the group $T$. 
For instance we identify a mapping class defined by an element $x$ of $T$ 
from the previous section with the element of $PPSL(2,\Z)$ that has the same 
action as $x$ on the triangulation of $H^2$ in which boundary circles 
of $D$ are crushed onto the vertices of the Farey triangulation.  
Using this identification between $T$ and $PPSL(2,\Z)$ we can state:   

\vspace{0.2cm} 
\noindent 

\begin{lemma}
The map $\Phi$ is  the unique anti-isomorphism between $T$ and $PT$
determined by the formulas: 
\[ \Phi(\alpha)=F, \;\; \Phi(\beta)=R\]
where $\alpha,\beta$ are the generators of $T$ from the previous section. 
\end{lemma}
\begin{proof}
The local moves can act far way by means of conjugacies.
One associates to the local move $x$ the element $\Phi(x)\in T$. 
If we want to compute the action of  $\Phi(\gamma\cdot x)$  
we compute first the action of $\Phi(x)$ and then we have to 
act by some transformation $\Lambda$ 
which has the same effect as $\Phi(\gamma)$
had on the initial triangulation. But the triangulation has been 
changed by means of $\Phi(x)$.  This means that the 
transformation  $\Lambda$ is therefore equal to 
$\Phi(x)\Phi(\gamma)\Phi(x^{-1})$. 
This implies that 
$\Phi(\gamma\cdot x)=\Phi(x)\Phi(\gamma)\Phi(x^{-1}) \Phi(x)=
\Phi(x)\Phi(\gamma)$. 
\end{proof}
\begin{remark}
This correspondence will be  essential  below. 
It enables us to express arbitrary flips on a triangulation in
terms of the local moves $F$ and $R$. 
Since the moves are local, small 
words will lead to small differences in the triangulations.  
Eventually, we can translate 
(by means of the canonical anti-isomorphism
which reverse the order of letters in a word) any word in the  
generators $R$ and $F$ into an element of the group $T$, viewed 
as a word in the standard generators $\alpha$ and $\beta$.  
It is more difficult 
to understand the properties of a combing in terms of the action 
of $\alpha$ and $\beta$ on triangulations since the action is not 
local, and thus a short word might have a quite large effect on the 
combinatorics of the triangulation. 
\end{remark}

\subsection{Mosher's normal form for elements of $T$ on infinitely 
many flips}
Mosher proved that mapping class groups  of finite surfaces are   
automatic (\cite{Mo}). One might expect then that mapping class groups 
of {\em infinite surfaces} share also some properties closed to the 
automaticity, but suitably weakened by the infiniteness 
assumption. 

\vspace{0.2cm}
\noindent 
The aim of this section is to define a
first  natural combing for $T$ derived from Mosher's normal form. 
Unfortunately, this combing is unbounded. We will show next 
that it can be modified so that the new combing is asynchronously bounded.

\vspace{0.2cm}
\noindent 
Mosher's proof of automaticity consists of embedding 
the mapping class group in the 
corresponding Ptolemy groupoid and derive normal forms (leading to 
combings) for the latter. The way to derive normal forms is however 
valid for all  kind  of surfaces, without restriction of their  
-- possibly infinite -- topology. The alphabet used by the automatic 
structure is based on the set of combinatorial types of flips. 
The only point where the finiteness was used by Mosher is when one 
observes that the number of different combinatorial flips on a 
triangulated surface (with fixed number of vertices) 
is finite, provided that the surface is finite. 
Thus, the same proof does not apply to the case of $T$, since 
there are infinitely many combinatorially distinct  flips. 
Nevertheless, we  already remarked that we can express 
an arbitrary flip on the infinite triangulation 
as a composition of the two elements $F$ and $R$. 
This observation will enable us to {\em rewrite} the Mosher   
combing in the Ptolemy group (which uses all flips as generators) 
as a new combing 
which uses only  the two generators $F$ and $R$. 
We will call it the {\em Mosher-type combing} of $T$. 
Recall that it is  
equivalent to a combing in our favorite 
generators $\alpha$ and $\beta$.

\vspace{0.2cm}
\noindent 
Let us recall the normal forms for elements of $T$, in terms 
of flips. Choose a base triangulation $\tau_B$, fixed once for all, 
for instance the Farey triangulation. Choose a total ordering on 
the edges of the triangulation $\tau_B$, say $e_1,e_2, e_3,\ldots$, 
so that $e_1=e$ is the d.o.e., and choose an  arbitrary orientation 
on each edge. The combing might depend on the particular choice we made. 
Given an element $\zeta\in PT$ we represent it as the 
couple of labeled triangulations $(\tau_B,\tau)$. 
These two triangulation are identical outside some finite 
polygon, on which the restriction of the two triangulations 
are different. Let us denote by $\tau_B^f$ and $\tau^f$ the 
restrictions of the two triangulations to the finite polygon. 
Then $\tau^f$ is obtained from $\tau_B^f$ by removing several 
disjoint edges, 
say $g_1=e_{j_1} < g_2=e_{j_2} <\ldots <g_k=e_{j_k}$,   
from $\tau_B$ and replacing them by another $k$  disjoint ideal 
arcs (which do not belong to $\tau_B$) having 
the same set of endpoints.

\vspace{0.2cm}
\noindent As it is well-known  any two  triangulations of 
a polygon could be obtained one from another by means of 
several flips. Moreover, if the  polygon had $n$ vertices
and thus it is partitioned into $(n-2)$ triangles then the
minimal number of flips needed to transform one triangulation
into another one is $2n-10$, and this estimation is 
sharp for large $n$, as it was proved by Sleator, Tarjan and 
Thurston (see \cite{STT}).  Notice that in our case 
the size of the 
polygon is not {\em  a priori} bounded since there exist elements 
of $T$ whose support might be arbitrarily large. 

\vspace{0.2cm}
\noindent 
Further we construct a series of flips followed by a relabeling 
move, the sequence being uniquely determined by the given 
element $(\tau_B,\tau)$. 
The edges $g_1, g_2, \ldots, g_k$ are called the uncombed 
edges (in this order). Consider the first one, namely $g_1$. 
We define a {\em  prong} of a triangulation 
to be a  germ of the angle determined by two 
edges of the triangulation which are incident, thus having 
a common vertex. 
We say that an oriented  ideal arc $g$  {\em belongs to some prong} of 
the triangulation $\tau^f$ if the start point of $g$  is the 
vertex of the prong and  $g$ is locally contained in the prong.

\begin{center}
\includegraphics{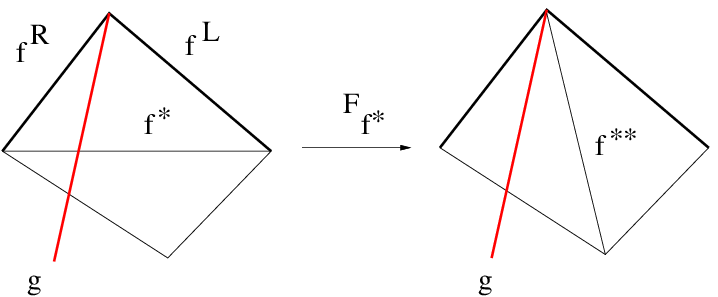}
\end{center}

\vspace{0.2cm}
\noindent 
Start with $g=g_1$ and recall that $g$ has been given an orientation. 
Pick up the unique prong of $\tau^f$ to which $g$ belongs. 
This prong is determined by the two edges $f^{L}$ and $f^{R}$, 
sitting on the left and respectively on the right of $g$. 
Therefore there exists another edge - say $f^*$ - 
of $\tau^f$ joining the two endpoints of 
$f^{L}$ and $f^{R}$, other than the vertex of the prong and 
forming a triangle $T$ with the prong edges.  
Further there exists another triangle $T^*$ of $\tau^f$ which shares 
the edge $f^*$ with $T$, but whose interior is disjoint from $T$, and 
thus it lies in the opposite half-plane determined by $f^*$. 
It is clear that $g$ intersects $f^*$ in one point. 

\vspace{0.2cm}
\noindent 
The first step in combing is to use the flip on the edge 
$f^*$, which replaces $f^*$ by the other diagonal $f^{**}$ 
in the quadrilateral $T\cup T^*$. 
If $g$ was precisely the diagonal  $f^{**}$ of $T\cup T^*$ then 
we succeeded in combing it, since it will belong to the 
new triangulation $\tau_1$, obtained by flipping. 
We will restart our procedure for $g=g_2$ and 
so on. Otherwise, it means that $g$ is still uncombed in the 
new triangulation. The  former prong determined by the arcs 
$f^{L}$ and $f^{R}$ is now split into the union of two 
prongs because we added one more arc, namely $f^{**}$ which shares 
the same vertex.  Then $g$ will belong to precisely one of the 
two new prongs, either to that determined by $f^{L}$ and $f^{**}$, or 
to that determined by $f^R$ and $f^{**}$. 
We change the notations for the arcs of the new 
 prong to which $g$ belongs so that the edge of the left is 
$f_2^{L}$ and the edge of the right is $f_2^{R}$. 
Then we restart the algorithm used above for $f_1=f$. Namely, 
consider the triangle $T_2$ determined by the two edges $f_2^L$ and 
$f_2^R$ and the edge $f_2^*$ connecting their endpoints, and next the 
opposite triangle $T_2^*$. Use the flip on the edge $f_2^*$, and 
continue this way. 

\vspace{0.2cm}
\noindent 
The lemma {\em Combing terminates} section 2.5 of \cite{Mo} 
tells us that after finitely many steps we obtain a triangulation 
for which $g_1$ is combed. We continue then by using the same 
procedure in combing $g_2$ and then $g_3$ and so on, until 
all $g_j$ (for $j\leq k$) are combed. 
At the end we need to relabel the d.o.e. in order to bring it 
to the d.o.e. of $\tau_B$. The  sequence of flips 
and  corresponding triangulations 
$\tau\to \tau_1\to\tau_2\to\cdots \tau_N\to\tau_B\to \tau_B$
(the last move being a relabelling) is  called the 
Mosher normal form of the element $\zeta$ of $T$.

\vspace{0.2cm}
\noindent 
As already mentioned before this normal form is not convenient 
for us  as it states, since there are infinitely 
many distinct  combinatorial flips.  We can 
overcome this difficulty by  translating in  the simplest possible 
way Mosher's normal form into a word in $F$ and $R$. 
In this case the flip $F$ cannot be applied but on the d.o.e. 
We assume that the d.o.e. $e_{\tau}$ of $\tau$ belongs to the 
polygon associated to $\tau_f$ and the same for $\tau_B$. 
This can be realized by enlarging the size of the support.  

\vspace{0.2cm}
\noindent 
In order to apply a flip to the triangulation $\tau_f$ we notice that 
we have first to move the d.o.e. from its initial position 
onto the edge which we want to be flipped. 
Let us assume for the moment that it is always possible to do this  
in a canonical way. Let $\tau$ be an arbitrary unlabeled triangulation
(finite or not) and $e, f$ two oriented edges. We define the 
{\em transfer} $T_{e,f}(\tau)$ as being the (unique) element 
of $T$ which sends the labeled triangulation $(\tau,e)$ 
into the labeled triangulation $(\tau,f)$. 

\vspace{0.2cm}
\noindent 
The normal form obtained above for $\zeta$ can be read now 
in the following way: 
\begin{enumerate}
\item Locate the first edge to be flipped, namely  
$f_1^*$, of $\tau^f$. 
Use the transfer $T_{e_{\tau}, f_1^*}(\tau^f)$ in order 
to move the d.o.e. from $e_{\tau}$ to $f_1^*$. 
\item Use the flip $F$, which will be located at $f_1^*$ and thus 
it will act exactly as the flip considered in the 
Mosher normal form.
\item The new d.o.e. is the image $F(f_1^*)$ of the former d.o.e.
with the d.o.e. orientation induced by the flip. 
Locate the new edge to be flipped, say $f_2^*$. Use the transfer 
$T_{f_1^*, f_2^*}(F_{f_1^*}(\tau^f))$. 
\item Continue until all uncombed edges are combed. 
\item If all edges were combed, then use eventually the transfer 
$T_{f_N^*, e}(F_{f_N^*}(\cdots (F_{f_1^*}(\tau^f))\cdots )$
in order to bring the d.o.e. at its right place. 
\end{enumerate}

\subsection{Writing Mosher's normal form as two-generator words}
We will explain now  how any transfer can be written 
{\em canonically} as a word in  the two letters $F$ and $R$, corresponding 
to the  respective generators of $T$. This procedure will be called 
then the {\em translation} of Mosher's normal form.  In fact, 
the transfer moves preserve the combinatorics of the triangulation, 
and thus they can be identified (by means of the anti-isomorphism
$\Phi^{-1}$ encountered above)  with automorphisms of the dual tree. 
Using this identification each transfer $T_{e,f}(\tau)$ 
corresponds to the element of the modular group 
$PSL(2,\Z)$ whose action on the binary tree 
sends the edge dual to $e$ to the edge dual to $f$. 
Now, it is well-known that $PSL(2,\Z)$ has an automatic 
structure and thus any element can be given a normal 
form in the standard generators. However, one can do this 
in an explicit elementary way. 
In fact the subgroup $PSL(2,\Z)$ of $T$ is actually
the (sub)group generated by the elements 
$\alpha^2$ and $\beta$. Furthermore, since we have 
$PSL(2,\Z)=\Z/2\Z *\Z/3\Z$, it follows that any element 
of $PSL(2,\Z)$ 
can be uniquely written as a word 
\[ \beta^{\epsilon_0}\alpha^2\beta^{\epsilon_1}\alpha^2\cdots 
\alpha^2\beta^{\epsilon_{m+1}}, \;\; {\rm 
where } \,\; \epsilon_0,\epsilon_{m+1}\in\{0,1,2\} \;\,{\rm and }\;\; 
\epsilon_1,\epsilon_2,\ldots, \epsilon_m\in\{1,2\}\]
 
\vspace{0.2cm}
\noindent 
Moreover, the factors $\epsilon_j$ can be effectively computed. 
There exist a completely analogous description in terms of the 
triangulations and the moves $F^2$ and $R$. We will explain this 
in the dual setting (since anyway the final result can be 
easily recorded as a word in $\alpha$ and $\beta)$. 
The ideal arcs $e$ and $f$ correspond to the edges of the dual binary 
tree ${\rm Tree}(\tau)$  of the triangulation $\tau$. 
Recall that all edges of $\tau$ have been given an orientation. 
This induces a co-orientation of the edges of  ${\rm Tree}(\tau)$, 
namely a unit vector orthogonal to each 
edge of  ${\rm Tree}(\tau)$.
One might choose a natural co-orientation for  ${\rm Tree}(\tau)$
by asking it to turn clockwisely in the standard planar 
embedding of the binary tree, in which case the formulas 
are simpler. However, it would be preferable to do the computations 
in the general situation.  
Moreover,  ${\rm Tree}(\tau)$ is a rooted tree, whose root is the 
vertex of ${\rm Tree}(\tau)$ sitting in the right of the edge 
dual to the d.o.e. $e_{\tau}=e$.  In order to fix it we  used 
the co-orientation. 
Moreover, for any edge $g$ of the tree and chosen vertex $v$ of $g$ 
it makes sense to speak about the two other edges incident to $g$ at 
$v$, which are: one at the left of $g$ and the other one at 
the right of $g$. 
This follows from the natural circular order around each 
vertex inherited from the embedding of the dual tree in 
the plane. 
If $e$ and $f$ are  -- not necessarily distinct --  
edges incident at some vertex,  we set: 
\[\epsilon(e,f)=\left\{\begin{array}{ll}
0 & {\rm if} \; e=f \\
1 & {\rm if} \; f  \; {\rm is\; on\; the\; left\; of} \;  e \\
2 &  {\rm otherwise}
\end{array}
\right.
\]
Furthermore, if we identify $\alpha^2$ and $\beta$ with elements 
of $PSL(2,\Z)$ which act as  planar tree automorphisms, 
it makes sense to look at the 
image of the co-orientation of an edge $e$ by means of the element 
$\xi\in PSL(2,\Z)$. For instance $\alpha^2$ reverts the orientation 
of the d.o.e. 

\begin{center}
\includegraphics{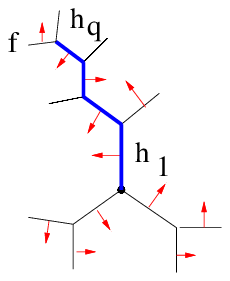}
\end{center}

\vspace{0.2cm}
\noindent 
Let then $\gamma=(h_1,h_2,\ldots,h_q)$ be the unique geodesic in 
${\rm Tree}(\tau)$ which joins  the root (which is an endpoint 
of $e$) to $f$. It might happen that  either $h_1=e$ or $h_1\neq e$, 
but  in any case we have always $h_q\neq f$. 
We claim now that: 
\[ \Phi^{-1}(T_{e,f}(\tau))=
\beta^{\epsilon_0}\alpha^2\beta^{\epsilon_1}\alpha^2\cdots 
\beta^{\epsilon_q}\alpha^{2\delta}\]
where 
\[ \epsilon_0=\epsilon(e, h_1), \;
\epsilon_1=\epsilon(h_1, h_2), \ldots
\epsilon_{q}= \epsilon(h_{q}, f)\]
and 
\[\delta=\left\{\begin{array}{ll}
0 & {\rm if} \;  f\; \mbox{has the co-orientation induced from } e \\
1 &  {\rm otherwise}
\end{array}
\right.
\]
Remark that the intermediary co-orientations of $h_j$ do not 
influence the normal form. 

\begin{remark}
If the triangulated  subpolygon  of the Farey triangulation 
is connected then it is  actually a convex polygon in the plane. 
In particular, if $e,f$ are edges of some of its triangulations and 
$h_j$ are edges dual to a geodesic joining the two dual edges,  then 
all $h_j$ are contained in the respective  subpolygon.
Thus, in the process of of combing uncombed edges we can realize 
all flips and transfers within the given, fixed polygon. 
\end{remark}

\begin{definition}
The Mosher-type combing (or normal form) for elements of 
$T$ is defined as follows. For each 
$\zeta=(\tau_B,\tau)\in T$ we choose the restricted triangulations 
$\tau_B^f, \tau^f$ so that the associated polygon is the 
smallest connected polygon containing all uncombed arcs of 
$\tau$ and both d.o.e.'s. Notice that this is uniquely determined. 
The normal form of $\zeta\in T$  is then the  
sequence 
\[\sigma^M_{\zeta}= 
[T_{e_{\tau},f_1^*}(\tau^f),\; F,\; T_{f_1^*,f_2^*}(F_{f_1^*}(\tau)),\; F,\; 
\ldots\;,\;F,\; T_{f_N^*, e}(F_{f_N^*}(\cdots (F_{f_1^*}(\tau^f))\cdots )]\]
in which each transfer is translated as a canonical word 
in $R$ and $F$. We might eventually use $\Phi$ in order 
to uncover the word in $\alpha$ and $\beta$. 
\end{definition}

\vspace{0.2cm}
\noindent 
Recall that Mosher's combing of the mapping class group 
is asynchronously bounded 
in the case when the surface is finite. 
This follows from the fact that, given two elements 
$\zeta,\widetilde{\zeta}$ at distance one in the Cayley graph 
of the  respective Ptolemy groupoid, then one can write  
\[\sigma^M_{\zeta}=w_0z_1w_1z_2w_2\cdots z_pw_p, \;\;
\sigma^M_{\widetilde{\zeta}}=\widetilde{w_0}z_1\widetilde{w_1}z_2
\widetilde{w_2}\cdots z_p\widetilde{w_p}\]
where $z_i,w_i$ are words in the generators, such that: 
\begin{enumerate}
\item First, the  size of the subwords on which the combings 
don't agree are uniformly bounded:
\[ |w_i|, |\widetilde{w_i}| \leq K, \;\; {\rm for\; all } 
\;\; i\leq p \]
\item Second, for all $m$,  
the distance between the corresponding {\em prefix elements} 
which are represented by the prefix words
$\sigma^M_{\zeta}(m)=w_0z_1w_1z_2w_2\cdots z_m$
and $\sigma^M_{\widetilde{\zeta}}(m)=\widetilde{w_0}z_1\widetilde{w_1}z_2
\widetilde{w_2}\cdots z_m$ is {\em uniformly bounded}, 
when  these elements are considered in the Cayley graph 
of the Ptolemy groupoid. 
Actually, the two prefix elements are always at distance one, 
because they differ by precisely one flip. 
\end{enumerate}
\noindent Let us analyze what happens in the  infinite case, 
when the Ptolemy group is $T$. We rewrite the combing 
in  the two-generator free group and thus we have to 
rewrite also the flip relying  the two prefix elements 
above.  
Since there are flips which are arbitrarily far away, we will need 
arbitrary long words in $R$ and $F$, 
and thus the Mosher-type combing is not 
asynchronously bounded. 
Notice however, that the first part of the assertion above is still true 
in this case. 
This follows from the fact that $\zeta, \widetilde{\zeta}$ are at 
distance one in the two-generator Cayley graph if they differ by a 
fundamental flip or by a rotation move $R$. In the second case 
the combing paths are the same except for the first few moves 
which realize the first transfer.

\subsection{Modifying the Mosher-type combing in order to get asynchronous boundedness}
We turn  back to the original description of the  
Mosher combing of $T$ in terms of the infinitely many 
flips of the triangulation. Any element
$\zeta=(\tau_B,\tau)$ of $T$ was brought to its
normal form 
\[ F_{e_1},F_{e_2},\ldots, F_{e_N}, P\]
where $F_{e_i}$ are flips and $P$ is the last relabelling move. 
We will  first define the new combing -- denoted
$\sigma$  -- in the usual generators $R$ and $F$ for 
each flip $F_{e}$ and then concatenate the combings  
according to the pattern of the Mosher normal form above. 
Recall that the Mosher-type combing  defined  above 
for the flip $F_f$ was the simplest possible (actually geodesic): 
we used the transfer of the d.o.e. $e$ to $f$,
further the fundamental flip $F$,  
and then we transfered the d.o.e. $F(f)$ back to its initial position
$e$. This time we will be more careful about the way we will 
achieve the flip $F_f$.  

\vspace{0.2cm}
\noindent 
Consider  the element $\widetilde\zeta$ of $T$ at distance 
$1$ from $\zeta$. Let us assume that 
$\widetilde\zeta=F\zeta$. Then 
$\widetilde\zeta=(\tau_B,\widetilde\tau)$, where 
$\widetilde{\tau}=F\tau$ is the triangulation with the 
d.o.e. $F(e)=e^*$. 

\vspace{0.2cm}
\noindent 
The failure of the boundedness for the 
Mosher-type combing above is a consequence of the fact 
that the distance between $T_{e,f}(\tau)$ and 
$T_{e^*,f}(\widetilde\tau)$ grows linearly with the distance 
$d(e,f)$ between the edges $e$ and $f$. 

\vspace{0.2cm}
\noindent 
We want to define  a path $P(\tau, F_f(\tau))$
joining $\tau$ to  $F_f(\tau)$ such that  the asynchronous distance 
between  the paths $P(\tau, F_f(\tau))$ and 
$P(\widetilde\tau, F_f(\widetilde\tau))$ is uniformly bounded, 
independently on the position of $f$. 

\vspace{0.2cm}
\noindent 
Let us denote by $Q_f$ the quadrilateral determined by the edge $f$, 
which has $f$ and $f^*$ 
as diagonals. We have two distinct cases to analyze: either  
$f$ is disjoint from $Q_e$ or $f$ belongs to $Q_e$. 

\vspace{0.2cm}
{\bf (1.)} Assume that $f$ is disjoint from $Q=Q_e$. 
Consider the chain of triangles joining $e$ to $f$, which is dual 
of the geodesic joining $e$ to $f$ in the dual tree. Denote their 
union by $Z(e,f)\subset \tau$.

\begin{lemma}\label{combflip} 
There exists a sequence of labeled triangulations 
$\tau_n$ with d.o.e. $e_n$ and a sequence of 
polygons $Z_n\subset Z(e,f)$ 
such that
\begin{enumerate}
\item  For all $n$ we have  $e_n, e_{n+1}\in Z_n$;
\item We have $Q \subset Z_n$ for all $n$;   
\item The number of vertices 
of $Z_n$ is uniformly bounded by $K$. We will see that $K=7$ suffices;  
\item For $n$ large enough $Z_n= Q\cup  W_n$, where 
$W_n$ and $Q_f$ have a common edge.  
\end{enumerate}  
\end{lemma}
\begin{proof}
Consider thus the triangulated polygon $A=Z(e,f)-Q\cup Q_f$, which can be 
seen as a chain of triangles joining the edges 
$e'\in Q$ (adjacent to $e$) and $f'\in Q_f$ (adjacent to $f$).
This chain of triangles is dual to a geodesic in the binary dual 
tree and thus it contains the minimal number of possible triangles. 
This polygon is embedded in the plane and it makes sense to speak 
about the left vertex  $a_1$ of $e'$ and the right vertex $b_1$ of $e'$. 
If we start traveling  in the clockwise direction along 
the boundary of $A$  and starting at 
$a_1$ then we will encounter, in  this order, the vertices 
$a_2,...,a_s$, the last one  being a vertex of $f'$. Further if we travel 
in the counterclockwise direction from $b_1$ we will encounter  
the vertices $b_2,b_3,\ldots,b_t$, the last point being the other vertex
of $f'$. If $s,t\leq 2$ then there is nothing to prove. Assume that 
$t\geq 3$, the other case being symmetric. 
Consider the points $b_1,b_2,b_3$ and the smallest triangulated 
subpolygon $B\subset A$ containing these three vertices. It is understood 
that the triangulation of $B$ is the restriction of that from $A$ and 
the adjective smallest means that it contains the minimum number of triangles. 
Notice that edges of the triangulation of $A$ cannot join two 
vertices on the same side, since otherwise the subpolygon 
they would determine (part of the boundary and this edge) could be 
removed from $A$, thus contradicting the minimality of the chain $A$. 

\vspace{0.2cm}
\noindent 
Assume that $B$ has the vertices $a_1,a_2,...,a_m$ from the left 
side and $b_1,b_2,b_3$ from the right side.  
 
\vspace{0.2cm}
\noindent 
If $m\geq 3$ then there exists a subpolygon $W_1$ containing $a_1b_1$
of height at most 3. This means a subpolygon containing 
three  consecutive points on the same side, say $b_1,b_2,b_3$,  
and at most two vertices consecutive vertices on the other side
say $a_1,a_2$.    
If $m\leq 2$ this is obvious. 
Suppose now that $m\geq 3$. 

\vspace{0.2cm}
\noindent 
Since $B$ has $m+3$ vertices one needs $m$ diagonals  
in order to triangulate $B$. Let us denote by $m_j$ the number of 
diagonals having $b_j$ as endpoint. We have 
then $m_1+m_2+m_3=m$. If $m_1\geq 2$ then the diagonals exiting 
$b_1$ should arrive at $a_2,a_3$ and thus the quadrilateral 
$a_1a_2a_3b_1$ has the claimed property. 
Suppose $m_1=1$; if  $m_2=1$  then the diagonals are 
$a_2b_1, a_2b_2,a_3b_2$ and the subpolygon is $a_1a_2a_3b_2b_1$. 
If $m_2=1$ then $m_3\geq 1$ and thus there exists the 
diagonal $a_2b_3$ and thus the subpolygon 
$a_1a_2b_3b_2b_1$ verifies the claim. 

\begin{center}
\includegraphics{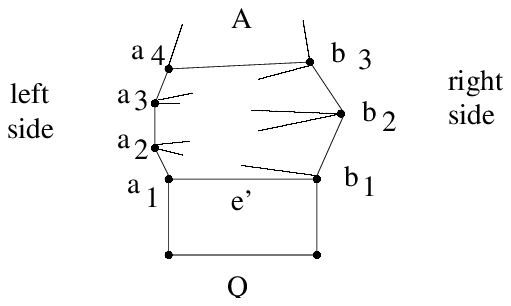}
\end{center}

\vspace{0.2cm}
\noindent 
Consider next the polygon $Z_1=Q\cup W_1$. Remark that 
$Z_1$ has  at most 7 vertices. We can use flips 
and transfers inside $Z_1$ (notice that the d.o.e. lays 
within $Q\subset Z_1$) in order to change the triangulation 
so that the tree consecutive points on one side of $W_1$ 
form now a triangle. We suppose that the d.o.e. is brought 
back into its position. 

\begin{center}
\includegraphics{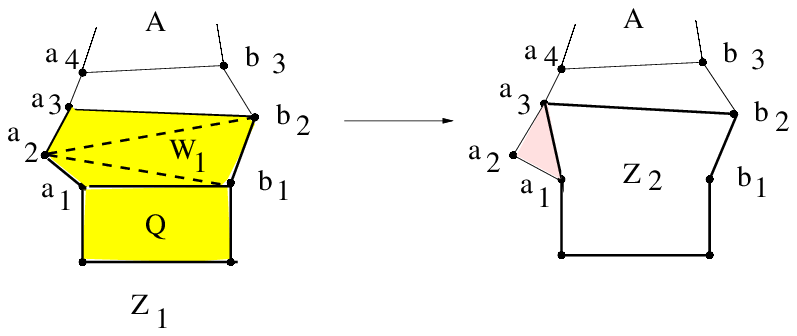}
\end{center}

\noindent 
Denote by $\tau_2$ the new triangulation. The chain of triangles 
$Z_2(e,f)$ which joins $e$ to $f$ within $\tau_2$ has one triangle less,  
because we can exclude the triangle $b_1b_2b_3$; 
actually one of the paths
$a_1,a_2,\ldots, a_s$ or $b_1,b_2,\ldots, b_s$
becomes one unit shorter. 

\vspace{0.2cm}
\noindent 
Continue the same procedure with the the polygon $A_2=
Z_2(e,f)-Q\cup Q_f$, and define inductively the polygons $W_2$,   
$Z_2=Q\cup W_2$ and so on, until we obtain a polygon 
$W_N$ that contains both $a_1b_1$ and $a_sb_t$.
This proves the claim. 
\end{proof}

\vspace{0.2cm}
\noindent 
Let us define now the combing of the flip $F_f$, as follows. 
The lemma \ref{combflip} shows that that there exist   
the subpolygons 
$Z_n$ constructed as above so that eventually 
$Z_n=Q\cup W_n$ where $W_n$ and $Q_f$ have a common edge. 
We will use the sequence $Z_j$ in order to join the d.o.e. $e$ to $f$.
However instead of going straight away from $e$ to $f$ along the 
shortest path we will consider the sequence $\tau_n$ of triangulations 
that makes  eventually the edges $e$ and $f$ to become close 
to each other. At each step, we consider only those 
flips or rotations which can be realized inside the respective $Z_j$ and 
that changes $\tau_j$ into $\tau_{j+1}$.  
Finally we obtain a triangulation containing $Z_n$. 
We will say that we modified the triangulation in order 
to make the transfer $T_{e,f}$  {\em short}. 
Now,  the polygon  $Z_n\cup Q_f$ contains both $Q$ and $Q_f$ and has 
at most 9 vertices. Therefore we  can realize  the 
flip $F_f$ within $Z_n\cup Q_f$, by using the transfer 
$T_{e,f}$ within the 9-vertices polygon, followed by  the 
fundamental flip $F$ and then  by the inverse transfer 
$T_{e,f}^{-1}$. The effect on the triangulation was precisely 
that of $F_f$. 

\vspace{0.2cm}
\noindent
Next, after the flip $F_f$ was performed, we use the 
the same procedure by means of the subpolygons $Z_n$ 
and the inverse sequence of associated triangulations 
in order to  move backwards  thru triangulations 
and finally reconstruct the original triangulation within 
the polygon $Z(e,f)$, without touching  anymore to $Q_f$. 
Eventually, we obtain a triangulation corresponding to 
$F_f\tau$.  

\begin{definition}
The path  $P(\tau, F_f\tau)$  consisting of  transformations 
(flips and rotations) between triangulations that join 
$\tau$ to $F_f\tau$ is the combing $\sigma_{F_f}$ of $F_f$. 
\end{definition}

\begin{remark}
There is some freedom in choosing the way to transform 
the triangulations within each $Z_n$. However, this does not 
influence the boundedness properties of the combing. 
We can easily find canonical representatives, since there are 
finitely many choices. 
\end{remark}

\begin{lemma}\label{bound}
Assume that $\widetilde\tau =F\tau$  differ by a flip 
$F$  and $f$ is an edge disjoint from $Q$. 
Then the combing paths $P(\tau, F_f\tau)$ and 
$P(\widetilde{\tau}, F_f(\widetilde{\tau}))$ 
stay at bounded asynchronous  distance in the  standard two generators 
Cayley graph of $T$. 
\end{lemma}
\begin{proof}
We have the sequences $Z_n, \widetilde{Z_n}$ as above. 
Then $Q_{e^*}=Q_{e}=Q$, and thus $Z_n\cap \widetilde{Z_n}\supset Q$
for all $n$. Further, for  $n$ large enough  we have 
$W_n=\widetilde{W_n}$ is a pentagon containing both an edge 
of $Q$ and an edge of $Q_f$.

\vspace{0.2cm}
\noindent 
Now, for any $n$ the polygon  $Z_n\cup \widetilde{Z_n}$
is connected and has at most 
$8$ vertices, since both lay in the same half space 
determined by the common edge $e'$  and thus have 
one more common vertex. This means that one can pass from 
$(\tau_n,e_n)$ to $(\widetilde{\tau_n}, \widetilde{e_n})$ 
by using only flips and  transfers taking place 
in  the finite polygon $Z_n\cup \widetilde{Z_n}$ and thus they are at 
bounded distance. For instance their distance is smaller than the 
diameter of the graph of transformations of a 8-vertices polygon 
(using $F$ and $R$) which is smaller than $30$. 
\end{proof}

\vspace{0.2cm}
{\bf (2.)} The second case to be considered is when $f$ is a nearby 
edge, namely  an edge of $Q$. 
\begin{definition}
If $f\in Q$ then $\sigma_{F_f}$ is the  
Mosher-type combing of $F_f$ in the two-generators, 
i.e. $\sigma^M_{F_f}$. 
\end{definition}

\vspace{0.2cm}
\noindent 
We are able now to define the combing of a general element of 
$T$. Let us consider $\zeta\in T$, which is written in Mosher's normal 
form in terms of arbitrary flips as a sequence 
$F_{f_1}, F_{f_2},...,F_{f_N}$ followed by a relabelling move $P$ 
bringing the d.o.e. at its place. 

\vspace{0.2cm}
\noindent 
We consider that the d.o.e. is at position $e$ and we 
replace each flip $F_{f_j}$ in the sequence above by its combing 
$P(\tau_j, F_{f_j}\tau_j)$, where $\tau_{j+1}=F_{f_j}\tau_j$. 
The d.o.e. remains at the same place $e$, except when  
$f_j=e$ in which case is transformed accordingly. 
At the end we obtain the 
triangulation $\tau_B$ with some place for the d.o.e. which is 
then transferred by means of the Mosher-type combing of 
the transfer onto its standard location.

\subsection{The combing of $T$  is asynchronously bounded}
The proof follows now the same lines as Mosher's proof. 
Consider $\zeta$ and $\widetilde\zeta$ two elements at distance one. 
We have either $\widetilde\zeta=F\zeta$ or $\widetilde\zeta=R\zeta$.

\begin{proposition}
The combing defined above for $T$ is  asynchronously bounded, 
and we can take the constant $K=30$. 
\end{proposition}
\begin{proof}
We have either $\widetilde\zeta=F\zeta$ or $\widetilde\zeta=R\zeta$.
It suffices to analyze the first case, the second case 
being simpler and resulting by the same argument. 
Let us consider the  Mosher normal forms 
\[ F_{f_1}, F_{f_2},\ldots, F_{f_n}, P ~ \; \mbox{\rm and respectively }~ \; 
\widetilde F_{f_1}, \widetilde F_{f_2},\ldots, \widetilde F_{f_n}, 
\widetilde P \]
From \cite{Mo} section 2.5 this combing is asynchronously bounded 
if we consider all  flips as generators, and moreover 
the combing sequences $F_{f_j}$ and $\widetilde{F_{f_j}}$ 
coincide at those positions corresponding 
to  flips outside the quadrilateral $Q$. Notice that our $F$ 
is located at the d.o.e. while  \cite{Mo} deals with 
the general case of the flip which can be outside the d.o.e. 
Therefore the idea  of the proof is  very simple: the points 
in the Mosher combing corresponding 
to the flips which are located at edges of $Q$ are at bounded 
distance from each other; this distance is measured by 
composing a few flips, which are themselves  flips
on edges uniformly closer to the d.o.e.  
Thus after transforming them into paths in the 
two generator Cayley graph these points will be only a bounded 
amount apart. The points corresponding to flips on edges which  
are far from the d.o.e.  could be  very far away 
in the Mosher-type combing, but these 
points come from identical sequences of flips and each flip has 
been combed now using the paths $P(\tau, F_f\tau)$.
The lemma \ref{combflip} shows that these points will remain also 
a finite amount apart. Eventually, we have to see what happens 
when using the relabelling moves $P, \widetilde{P}$. It suffices to 
observe that the d.o.e. of $\tau$ will remain always closed-by 
to the  d.o.e. of $\widetilde{\tau}$, and actually in the same 
quadrilateral. This means that the last transfer of d.o.e. 
leads to two normal forms which are very closed to each other. 
This will prove the claim. 

\vspace{0.2cm}
\noindent 
It suffices thus to see what happens with Mosher-type combing 
when we meet nearby edges to be flipped.

\vspace{0.2cm}
\noindent 
Let $f_1^*$ be the first edge to be flipped and 
$e=e_{\tau}$ be the d.o.e. We have to compare 
$T_{e,f_1^*}(\tau^f)$ and 
$T_{\beta(e),f_1^*}(\tau^f)$. An alternative way is to 
look at the dual tree. Recall that ideal arcs of the 
triangulation yield edges of the dual tree. 
We have to
compare the two geodesics $\gamma$ and $\widetilde{\gamma}$
which join the left endpoints of the edges $e$ and 
respectively $\beta(e)$  to some endpoint of $f_1^*$. 
But the left endpoints of $e$ and $\beta(e)$ coincide
and thus $\gamma=\widetilde{\gamma}$. 
Thus the transfer are given by identical words 
except possibly for the first three letters.

\vspace{0.2cm}
\noindent 
The next case is when $\widetilde{\zeta}=\alpha\zeta$. 
Then the differences between the normal forms can propagate 
to all other transfers and not just to the first one. 
Another difficulty is that the dual trees are different. 
We set 
${\mathcal T}={\rm Tree}(\tau)$ and 
${\mathcal T}_t$ for the dual tree after $t$ steps. 
We define a {\em step} to be the action of a block 
of several consecutive letters of the normal form. 
The precise control  on the size of blocks will be given below. 
We set also 
$\widetilde{{\mathcal T}}={\rm Tree}(\alpha\tau)$
and then $\widetilde{{\mathcal T}}_t$ for the 
dual tree after $t$ steps. The steps in the two cases are 
not necessary correlated.  Instead, we  would rather want 
a certain correlation between the trees 
${\mathcal T}_t$ and $\widetilde{{\mathcal T}}_t$, for any $t$. 

\vspace{0.2cm}
\noindent 
The trees ${\mathcal T}$ and $\widetilde{{\mathcal T}}$ are 
identical except for the  image of the 
support $\Sigma$ of the move $\alpha$, which is made of the 
edge $e$ and its four adjacent edges. We have 
$\widetilde{{\mathcal T}}={\mathcal T}-\Sigma \cup \alpha \Sigma$, 
where $\Sigma$ is replaced by its image  by $\alpha$ 
(i.e. a rotation of angle $\frac{\pi}{2}$).  
We would like to define the steps in such a way that 
any moment $t$ we have 
$\widetilde{{\mathcal T}}_t={\mathcal T}_t-\Sigma_t \cup 
\alpha \Sigma_t$, where $\Sigma_t$ is combinatorially 
isomorphic to $\Sigma$. We call $\Sigma_t$ the singular locus 
at step $t$ and denote $\widetilde{\Sigma_t}=\alpha\Sigma_t$. 
Moreover we have a natural combinatorial 
isomorphism between the two trees, outside their 
respective singular loci.  Let $e_t$ denotes the 
central edge of $\Sigma_t$ and $\widetilde{e_t}$ 
for $\widetilde{\Sigma_t}$.

\vspace{0.2cm}
\noindent In order to get control on the differences 
between the normal forms in the two cases we have to understand 
what happens if we have to use transfers or flips which 
touch the singular locus. 
In fact, any transfer between two edges lying in the same 
connected component of ${\mathcal T}_t-\Sigma_t$ has 
a counterpart as a transfer in $\widetilde{{\mathcal T}_t}-
\widetilde{\Sigma_t}$ given by the same word.

\vspace{0.2cm}
\noindent 
As we saw previously the transfer between two edges 
is determined by the geodesic joining the two edges. 
We have then to understand what happens when such a 
geodesic penetrates in the singular locus. We have also 
to consider the case when we encounter a  flip on an 
edge from the singular locus.   
There are a few cases to consider: 

\begin{enumerate}
\item {\em If the geodesic enters and exit the singular locus.}
Let $\gamma=(h_1,h_2,\ldots,h_q)$, and 
respectively $\widetilde{\gamma}=(\widetilde{h_1},
\widetilde{h_2},\ldots,\widetilde{h_q})$ be the two corresponding 
geodesics which join two edges $f$ and $g$ which are corresponding 
to each other and both lay outside the singular locus. 
It follows that $\gamma-\gamma\cap \Sigma_t=
\widetilde{\gamma}-\widetilde{\gamma}\cap \widetilde{\Sigma_t}$
and the only differences can be seen at the level of the 
singular loci. According to the formula for the transfer 
we can write then 
\[ T_{f,g}(\tau^f_t)=z_1w_1z_2, \;
T_{f,g}(\widetilde{\tau^f_t})=z_1\widetilde{w_1}z_2\]
where the words $w_1,\widetilde{w_1}$ record the 
transformations needed to transfer one edge to another  
within the singular locus. The longest such word is 
$\alpha^2\beta\alpha^2\beta^2\alpha^2$ and thus 
$|w_1|, |\widetilde{w_1}|\leq 9$. 
\item {\em If the geodesic enters the singular locus and does 
not exit, or a geodesic starts from the singular locus 
and exits.}  This means that we have a transfer from 
an edge outside the singular locus to an edge 
of the singular locus. 

\vspace{0.2cm}
\noindent 
If this transfer is the final operation and the normal 
form is achieved  for $\zeta$, then normal form 
of $\widetilde{\zeta}$ is obtained by flipping the edge 
$\widetilde{e_t}$. 

\vspace{0.2cm}
\noindent 
Otherwise we didn't reach yet the normal form in neither
of the two configurations.  Thus the transfer is followed 
by a flip on some edge in the singular locus. We have two 
subcases: 
\begin{enumerate}
\item {\em The flip acts on some  edge $f$  of the 
singular locus incident but different from to $e_t$ in $\tau^f$ 
and different from $\widetilde{e_t}$ in $\widetilde{\tau^f}$.} 
Recall that we flip an edge 
in order to comb an uncombed ideal arc  $g$ which belongs to 
one of the two prongs determined by that edge. However, 
the ideal arc $g$ to be combed should belong to the 
prong opposite to the edge $e_t$. In fact, if $g$ belonged
to the prong containing $e_t$ then $g$ would intersect
(in the other picture, that of $\widetilde{\tau^f}$) 
first the edge $\widetilde{e_t}$. Thus the first flip 
in the process of combing $g$ would be the flip on the 
edge $\widetilde{e_t}$, contradicting our assumptions.

The possible situations are drawn below. We use now (for
a better intuition) the picture on the triangulation rather 
than on the dual.  
\begin{center}
\includegraphics{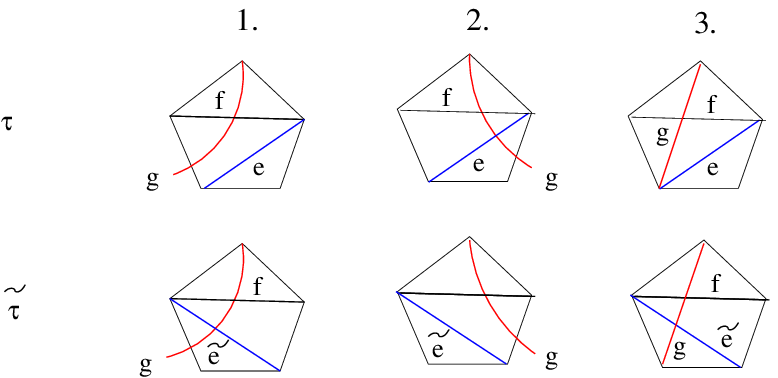}
\end{center} 
Let us analyze the first case, the third one being symmetric. 
The normal form reduction of $\tau^f$ takes the following form
and then it continues by combing the  ideal arc $g$ along 
the edge $h$. The d.o.e. is marked by a little square. 
We set further $e_{t+1}=e_t$ and $\Sigma_{t+1}=\Sigma_t$. 
\begin{center}
\includegraphics{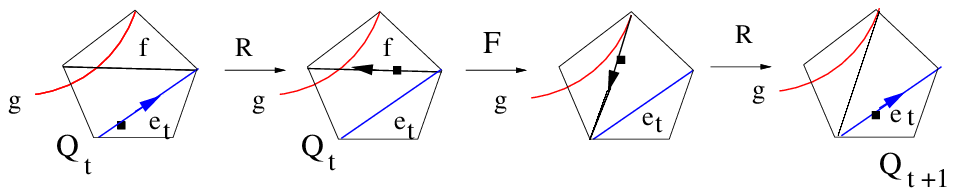}
\end{center} 
On the other hand the normal form reduction for 
$\widetilde{\tau^f}$ takes the form below and then continues 
by combing the arc $g$ along the edge $h$: 
\begin{center}
\includegraphics{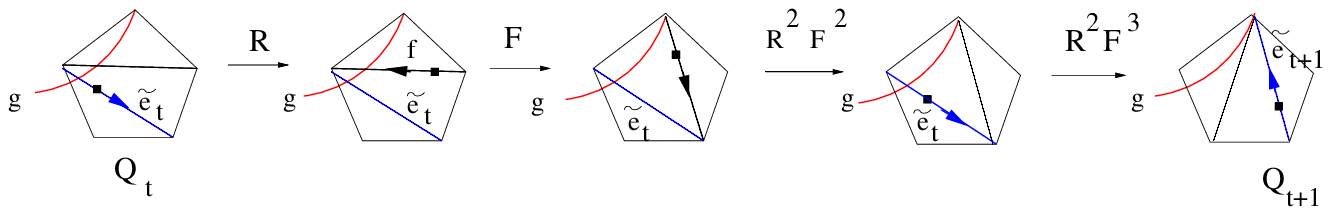}
\end{center}
Here one used the  fundamental flip $F_{\widetilde{e_t}}$ 
and this way the singular locus has been changed. Denote by 
$\widetilde{e_{t+1}}$ the  other edge in the pentagon 
as in the figure above and set $\widetilde{\Sigma_{t+1}}$ for the 
quadrilateral with diagonal $\widetilde{e_{t+1}}$. 
It is obvious that $\Sigma_{t+1}$ and 
$\widetilde{\Sigma_{t+1}}$ correspond to each other by means of 
a flip on the edge $e_t$ and they form the singular loci of 
the respective  couple of triangulations. 

In both situations the normal form reductions are identical 
for now on. This means that there are subwords 
$u$ of the combing $\sigma_{\zeta}$ and 
$\widetilde{u}$ of the combing of $\sigma_{\widetilde{\zeta}}$
so that, we can read off the strings which might be 
different from the picture above (recall that the 
local moves strings should be read in reverse order):  
\[ u=z_1w_1z_2, \; \widetilde{u}=z_1\widetilde{w_1}z_2\]
\[ w_1= RFR, \; \widetilde{w_1}=R^2F^3R^2F^3R\]
Thus the subwords $u$ and $\widetilde{u}$ are identical except 
for an extra string of length 4 in $\tilde{u}$.  

A similar computation shows that in the second case we have 
the previous transformation finish the combing of $g$, and 
thus we have to look at the next ideal arc to be combed. 

\item The flip is on the edge $e_t$. The picture are 
similar to those from above. We skip the details. 
\end{enumerate}
\end{enumerate}
This ends the proof of the proposition. 
\end{proof}
\begin{remark}
The distance between 
the words formed by the first $t$ letters in the combing 
of $\zeta$ and $\widetilde{\zeta}$ is bounded 
by a function linear in $t$. In fact each time that we are 
crossing the singular locus (by example in a transfer) 
the distance may have a jump by some $k\leq 9$, and the number 
of such crossing can grow linearly with the length of the 
word. 
\end{remark}

\section{Combing the braided Thompson group}

\subsection{Generators for $T^{*}$}
It is known (see \cite{FK2}) that the group $T^*$ is also 
generated by two elements 
that correspond to $\alpha$ and $\beta$ above. 

\vspace{0.2cm}
\noindent 
Specifically, we consider the following elements of $T^{*}$: 
\begin{itemize} 
\item The support of the element $\bps$ is the central hexagon.  
Further $\beta$ acts as the counterclockwise rotation of order  
three whose axis is vertical and which permutes cyclically  
the punctures.  

\begin{center} 
\includegraphics{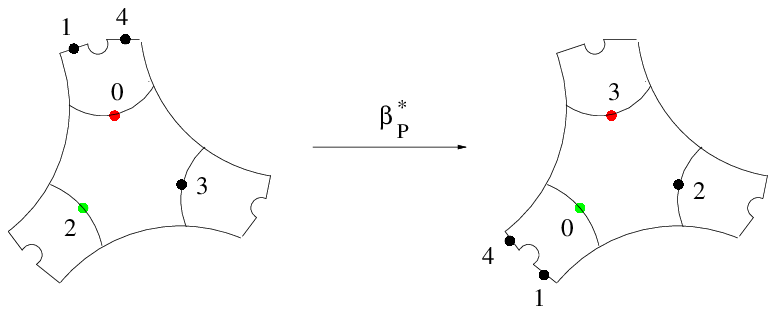}      
\end{center}
\item  The support of $\aps$  is the union of two adjacent hexagons, 
one of them being the support of $\bps$.  
Then $\aps$  
rotates counterclockwise the support 
of angle $\frac{\pi}{2}$, by  
keeping fixed the central puncture.  

\begin{center} 
\includegraphics{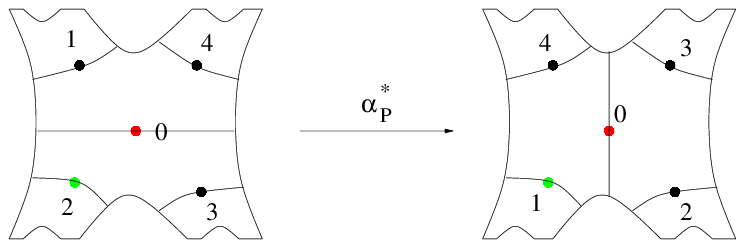}      
\end{center}
\end{itemize} 

It is proved in \cite{FK2} that   
$T^{*}$ is generated by $\aps$ and $\bps$.

\subsection{Normal forms for elements $T^*$}
The purpose of this section is to find a combing for $T^*$ 
using the generators $\aps,\bps$. The main novelty consists 
in using methods typical for mapping class groups that generalize 
first to $T$ and then to $T^*$. 
The main result of this section is 
\begin{theorem}
The group $T^*$ is asynchronously combable with departure function. 
\end{theorem}
\noindent From \cite{Ger} we obtain that 
\begin{corollary}
The group $T^*$ is ${\rm FP}_3$ and it has solvable word problem. 
\end{corollary}
\begin{remark}
It is claimed in \cite{Ger} that group like in the statement 
are  actually ${\rm FP}_{\infty}$, but the proof has not yet appeared 
in print. Another approach to the property ${\rm FP}_{\infty}$
is Farley's proof of the finiteness for braided picture groups. 
The group $T^*$ is a kind of picture group, where the role 
of permutations is now taken by the braid groups. 
Brin and Meier announced that this approach could lead to the proof 
of ${\rm FP}_3$ for the Brin group. 
\end{remark}

\subsection{The punctured Ptolemy groupoid  $T^*$}
In this section we will explain which are the modifications 
necessary for adapting the previous proof for $T$ to the 
case of the group $T^*$. 

\vspace{0.2cm}
\noindent The first observation is that we can view $T^*$ as 
a group of flip transformation on certain generalized triangulations 
of the punctured  hyperbolic plane, that will be called 
punctured triangulations or decompositions.  
Specifically, let us consider 
the Farey tesselation $\tau_B$  of the hyperbolic plane 
(in the disk model) and assume that we puncture each ideal arc
at its midpoint.  We will obtain now a triangulation whose 
edges are ideal arcs constrained to pass thru the punctures, as 
below:

\begin{center}
\includegraphics{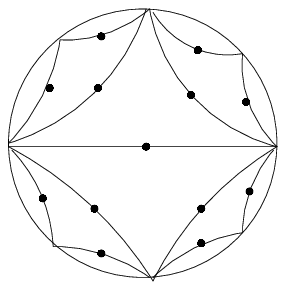}
\end{center}
We consider now that the punctures are fixed once for all.   
Then there is a set of moves which transform one such 
punctured triangulation into another one of the same type modelled 
on the transformations $F$ and $R$. We have the flip $F_{\gamma}$
on the punctured edge $\gamma$ and the rotation $R$ which changes 
the d.o.e. by moving it counterclockwise in the  (punctured) 
triangle sitting to its left:

\begin{center}
\includegraphics{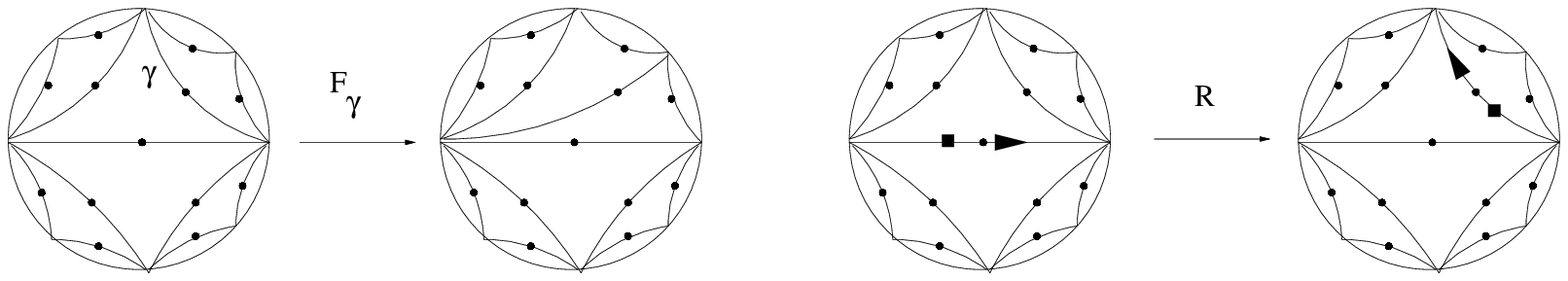}
\end{center}
Despite the similarities with the description of $T$ the fact that 
the punctured are fixed forces now  the ideal arcs to be distorted
and they cannot be realized anymore as geodesics in the 
hyperbolic plane.  For example, here is $(RF)^5$, where 
$F$ is the fundamental flip (on the d.o.e.): 

\begin{center}
\includegraphics{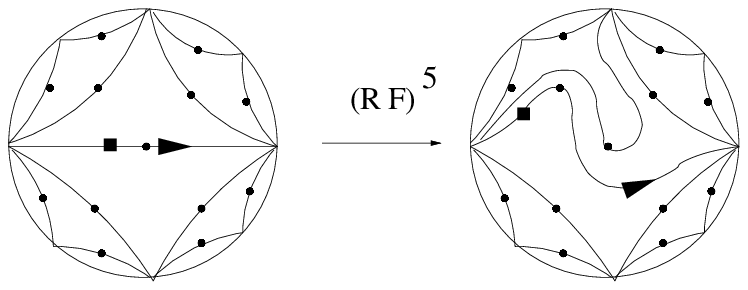}
\end{center}
We have then an immediate lemma:
\begin{lemma}
The punctured Ptolemy groupoid, which is the groupoid generated by  
flips on punctured triangulations is 
anti-isomorphic to the group $T^*$, by means of the 
anti-isomorphism $\Phi(\aps)=F$, $\Phi(\bps)=R$. 
\end{lemma}
\begin{proof}
Any flip can be realized as a composition 
of moves $F$ and $R$. 
\end{proof}

\begin{remark}
Each punctured arc of our triangulation  
splits into two half-arcs separated by a puncture. 
This way our punctured triangulations can be viewed as 
a hexagonal decomposition of the plane, each hexagon having 
three vertices at infinity and three more vertices among punctures. 
The decomposition could be refined to a triangulation by adding 
three extra edges in each hexagon, for instance the edges connecting 
pairwise the vertices at infinity. One can further consider the 
group generated by  all flips in the refined triangulation, which is 
the {\em Ptolemy groupoid of the punctured surface}. Notice that  
this group contains the 
{\em punctured Ptolemy groupoid $PT^*$}  defined above, as a proper 
subgroup. 
\end{remark}

\vspace{0.2cm}
\noindent The next task is to find normal 
forms in the punctured Ptolemy groupoid. Let us analyze 
what happens when using Mosher's  combing algorithm  
in the mapping class group of the punctured surface. 
First, a flipped arc should avoid all but one
punctures and thus it may not be represented by a 
geodesic in the hyperbolic plane or a  
straight segment in the flat plane.  

\vspace{0.2cm}
\noindent
This problem arose also in the case of the usual 
Ptolemy groupoid associated 
to an ideal triangulation of the punctured surface (see the 
remark above).  The solution given in that case is to  
consider only {\em tight} triangulations. Recall that two 
arcs are {\em tight} (with respect to each other) if they do not  
contain subarcs  bounding  a bigon i.e. an embedded  2-disk. 
Two triangulations are {\em tight} if  all their respective 
arcs are tight.  
Notice that we can pull triangulations tight and any flip can be 
realized as a flip between tight triangulations 
(see \cite{Mo}, section 2.5); thus the algorithm leading to  
normal forms works for the usual Ptolemy groupoid of 
the punctured surface.

\vspace{0.2cm}
\noindent 
However,  our decomposition is not a genuine triangulation 
of the punctured surface (but  rather a hexagon decomposition).  
In this respect,  the tightness of arcs  is not  concerning 
only the half-arcs going from one puncture to a point at infinity
(as it would be the case when dealing with the Ptolemy groupoid 
of the punctured surface), but  rather the entire arc. 
In fact, there exist triangulations having all their half-arcs tight 
although they are not tight. 
The flips which aimed at combing these arcs 
using  Mosher's algorithm
might not decrease the number of intersections points 
with the crossed arcs 
and thus the combing algorithm does not 
terminate. Here is a typical case of a uncombed arc for which 
the use of a flip move is not suitable: 

\begin{center}
\includegraphics{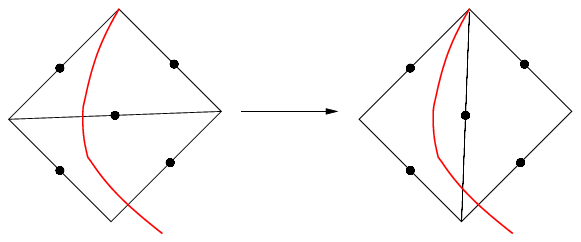}
\end{center}

\noindent 
In order to circumvent this difficulty we have to introduce  
some additional moves. 
We consider the braid twists $\sigma_{ef}$ 
which are elements of $T^*$  that are braiding (counterclockwise) 
the punctures sitting on the adjacent edges  $e$ and $f$.  
If $e$ is the d.o.e. 
then we can express $\sigma_{ef}$ as an explicit word 
in $F$ and $R$, which depends on the relative position of 
$f$ and $e$. However, in the next sections $\sigma_{ef}$ is simply  
a new letter in the alphabet. We want to find now an 
intermediary normal form of elements of $T^*$ using 
$F,R$ and the braid generators as well. 
Eventually, we will translate 
the obtained  normal forms as words in the two generators $R$ and 
$F$ alone. 

\subsection{Nonstraight arcs, conjugate punctures and untangling braids}

\paragraph{Straight arcs with respect to a  given triangulation.}
Let $(\tau_B^*,\tau^*)$ be an element of $T^*$. 
Thus $\tau^*$ coincides with the standard 
decomposition $\tau_B^*$ for all, but finitely many arcs. 
Consider an  ideal arc 
$\gamma$ which belongs to $\tau^*$, but not to $\tau_B^*$.
Our aim is to comb $\tau_B^*$ by means of flips and 
braid twists in order to transform it into a triangulation 
incorporating the arc $\gamma$. 
There are two situations. First, when $\gamma$ is isotopic
in  the disk $D$ (thus disregarding the punctures) to an arc of 
$\tau_B^*$ but there is no such isotopy 
which fixes the punctures (or, alternatively  
they are not isotopic in $D^*$). In this case we say that $\gamma$ 
is combed but it is not straight. This  
type of arcs should be {\em straightened}. 
The second possibility is that $\gamma$ is not isotopic 
in $D$ to an arc of $\tau_B^*$ and thus it has first to be combed and 
next to be straightened. 
We will give  below an algorithm which combs and straighten  
a given arc.

\vspace{0.2cm}
\noindent
Before to proceed, let us define properly what we mean by {\em 
straight edge}. We will work below with the flat planar model, 
but everything can be reformulated in the hyperbolic model. 
The triangulation  $\tau_B^*$ can be realized as a punctured 
triangulation of the disk, with vertices on the boundary circle 
(at infinity) called cusps and punctures in the interior of the disk. 
We assume that all edges are straight segments in the plane. 
Moreover, each edge is punctured at one point which is located 
at the intersection of the respective edge with the other diagonal 
of the unique quadrilateral to which the edge belongs. 
Let  now $\tau^*$ be an arbitrary triangulation which coincides 
with $\tau_B^*$ outside some finite polygon $P$. An edge of 
$\tau^*$ is {\em straight} if it is isotopic (in the punctured disk) 
to a straight segment (therefore keeping fixed the punctures).  
There is a similar notion which is defined using only combinatorial 
terms in the case of arcs inside a quadrilateral. 
Let $T$ be a triangle and $\gamma$ be a tight arc emerging from a vertex 
of it. We say that $\gamma|_{T}$ is {\em combinatorially 
straight} (or $\gamma$ is straight 
within $T$) if $\gamma$ intersects once more the boundary along the 
edge opposite to the vertex from which it emerges. Thus the different 
combinatorial models which might occur are those from below:

\begin{center}
\includegraphics{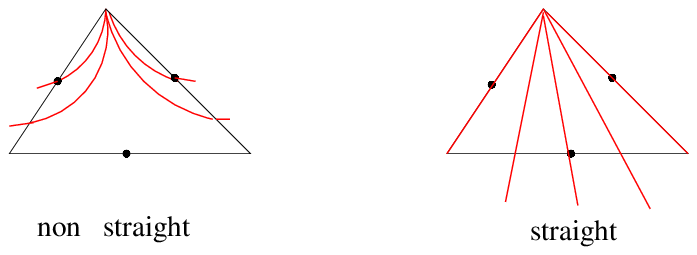}
\end{center}

\vspace{0.2cm}
\noindent
Let now consider now a quadrilateral $Q$ consisting of two 
adjacent triangles and $\gamma$ be a tight arc which 
emerges from one vertex 
$v$ of it. Let $q$ denotes the unique puncture inside $Q$.  
Let $v^*$ denote the vertex of $Q$ opposite to  $v$; then, the 
arc $vv^*$ punctured at $q$ splits $Q$ into two triangles $T^L$ and $T^R$. 
We say that  $\gamma|_{Q}$ is {\em combinatorially 
straight} (or $\gamma$ is straight 
within $Q$) if $\gamma$ is straight with respect to  both 
$T^L$ and $T^R$. This implies  that $\gamma\cap Q$ is contained either 
and $T^L$ or else in $T^R$. Typical examples of straight and non 
straight arcs are drawn below: 

\begin{center}
\includegraphics{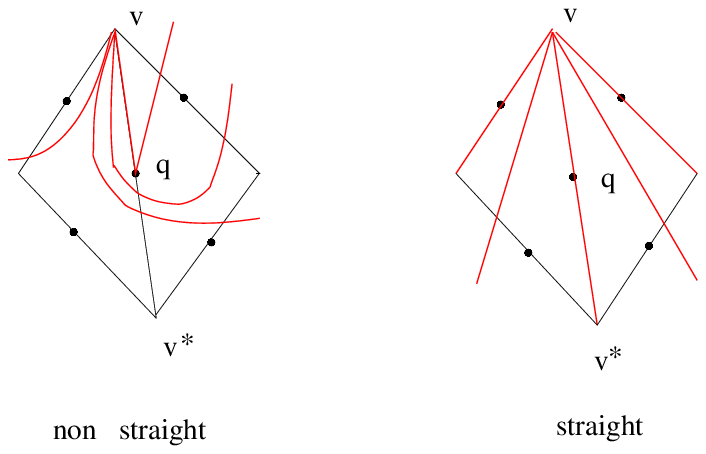}
\end{center}

\vspace{0.2cm}
\noindent Actually, as it can be seen an arc is combinatorially 
straight within a triangle or a quadrilateral if it can be isotoped 
by keeping its endpoints fixed to a line segment.

\vspace{0.2cm}
\noindent 
We consider from now on that all triangulations are isotoped so that 
their respective  half-arcs are tight.

\paragraph{Conjugate punctures along an arc.}
Recall that each edge $e$ of a triangulation  
has a puncture $p_e$ associated to it. If the triangulation is fixed then 
the puncture determines the edge.  

\vspace{0.2cm}
\noindent 
Consider now a tight arc $\gamma$ emerging from the vertex 
$v$ (opposite to the edge $e_0$) which crosses - in this order - 
the edges $e_1,e_2,...,e_m$ before ending in the vertex 
$v'$ (opposite to the edge $e_{m+1}$). 
Notice that the edge $e_j$ are not necessarily distinct. 
In order to determine completely the isotopy class of of the arc 
in the punctured plane one has to specify where sits the 
intersection point $\gamma \cap e_j$ with respect to the puncture
$p_{e_j}$. There are three possibilities, namely that the puncture 
be on the left side, on the right side or on the arc. We decide 
that the respective puncture is on the left side of $\gamma$ 
if this is so for an observer traveling along $\gamma$ 
in the direction given by the orientation of $\gamma$.  
We record this information by writing  down $p_{e_j}^L$. 
Similarly, when the puncture is on the right, we record this by 
writing $p_{e_j}^R$.  The superscript $L,R$ will be denoted   
$\epsilon(p_{e_j},\gamma)$  and called the {\em exponent
of the $j$-th puncture}. Notice that the same puncture might 
be encountered several times with different exponents. 

\vspace{0.2cm}
\noindent Extra caution is needed for the situation in which 
the arc pass thru the puncture $p=p_{e_j}$; 
this happens only once and only for one puncture, 
because the arcs we are interested in come from edges of 
punctured triangulations.   
We record this by adding an asterix $*$ to the letter $p_{e_j}$. 
Moreover, the exponent $\epsilon(p_{e_j},\gamma)$ can take now 
three values, namely from $\{L,R, 0\}$. 
Let $\gamma_{t}$ be an arbitrary small ${\mathcal C}^{\infty}$ 
perturbation of the arc $\gamma$ which is transversal to 
the arc $e_j$, is tight and avoids the puncture. 
Then we define $\epsilon(p_{e_j},\gamma)=\epsilon(p_{e_j},\gamma_t)$
for small $t$ if this is well-defined and 
 $\epsilon(p_{e_j},\gamma)=0$ otherwise. 
We can give more convenient ways to compute the exponent. 
\begin{lemma} 
\begin{enumerate}
\item If $e=vw$  and $p_e$ is the first  puncture encountered 
by $\gamma$ which emerges at $v$ then 
$\epsilon(p_e,\gamma)$ is $L$ if the frame $(p_ew,\gamma)$ 
is positively oriented and $R$ otherwise. Similarly if 
$p$ is the last puncture encountered by $\gamma$. 
\item 
Suppose that the local model of $\gamma$ is that of a 
local maximum at $p$ and thus the tangent vector 
$\dot{\gamma}_p$ points in the direction of $pw$ or else 
in the direction of $pv$. We consider that  
$\dot{\gamma}_p$ is a positive multiple of $pw$. 
Then $\epsilon(p_e,\gamma)=L$ if $\gamma$ lies locally 
on the left of the edge $vw$ (oriented as such) and 
$\epsilon(p_e,\gamma)= R$ otherwise. If 
$\dot{\gamma}_p$ is a positive multiple of $pv$ then 
the values  of the exponent are interchanged.  
\item 
Eventually in all other cases $\gamma$ crosses transversely 
the edge $e$ and cannot be reduced by isotopy to one 
of the previous two situations, and we set 
 $\epsilon(p_e,\gamma)=0$. 
\end{enumerate}
\end{lemma}
\begin{proof}
In the first two cases there are natural 
tight smooth perturbations giving the claimed values, while 
in the latter there both values could be reached by 
suitable perturbations:

\begin{center}
\includegraphics{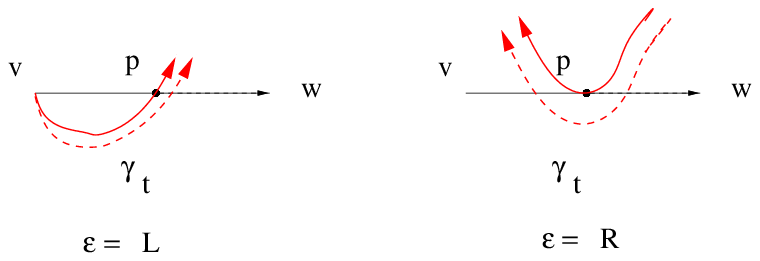}
 \end{center}
\end{proof}

\noindent If the exponent of a puncture is $0$ we say that it is 
an inert puncture. 
The arcs we are interested in come from edges of 
punctured triangulations and so they contain precisely one puncture. 
The word 
$w(\gamma)=v p_{e_1}^{\epsilon_1} p_{e_2}^{\epsilon_2}\cdots p_{e_m}^{\epsilon_m}v'$
where all $\epsilon_j\in\{L,R\}$ (excepting for one $j$ for which 
$\epsilon_j=0$) determines completely the isotopy class of the 
arc $\gamma$. Moreover, we will also say that $p(e_j)$ are the 
punctures that the arc $\gamma$ encounters. Furthermore, we 
define $p_k(\gamma)$ as the $k$-th puncture encountered by $\gamma$.

\begin{definition}
Let $p=p_i(\gamma)$ and $q=p_j(\gamma)$ be two punctures 
encountered by $\gamma$. We assume that none of them is  
an inert puncture. 
Then $p$ and $q$ are conjugate along $\gamma$ if the subword 
of $w(\gamma)$ starting at $p_{e_i}^{\epsilon_i}$ and ending at 
$p_{e_j}^{\epsilon_j}$ has the form 
$p_{e_i}^{\epsilon_i} p_{e_{i+1}}^{\epsilon_{i+1}}\cdots 
p_{e_j}^{\epsilon_j}$, where 
$\epsilon_i=\epsilon_{i+1}=\cdots =\epsilon_{j-1}\neq \epsilon_j\in
\{L,R\}$.  
In other words the punctures $p_i(\gamma),\ldots p_{j-1}(\gamma)$ 
stay on the same side with respect to $\gamma$ while the next 
puncture is the first one to stay on the opposite side. 
If one puncture, say $q$,  stays on $\gamma$ 
then we add extra condition, as follows.    
We ask that $q$ is the last puncture with exponent different 
from that of $p$ which  is encountered by $\gamma$.  
\end{definition}

\vspace{0.2cm}
\noindent Consider now two punctures $p=p_i(\gamma)$ and $q=p_j(\gamma)$
conjugated along $\gamma$. Notice that it might happen that the 
punctures $p_j$ are not distinct. 

\begin{lemma}
If $\gamma$ is tight then $w(\gamma)$ does not contain 
neither subwords of the form $p^{\epsilon}p^{\delta}$ with 
$\epsilon\neq \delta$, nor subwords of the form 
$p^{\epsilon}p^{\epsilon}p^{\epsilon}$. 
\end{lemma}
\begin{proof}
A subword of the form $p^{\epsilon}p^{\delta}$ corresponds 
to a subarc which is not tight and could be simplified 
by means of some isotopy. Further, an subarc corresponding 
to $p^{\epsilon}p^{\epsilon}p^{\epsilon}$ turns at least 
$2\pi$ around the puncture and its winding number with respect to 
$p$ is at least $2\pi$. 
However the winding number of the entire arc should be 
less than $\pi$. 

\begin{center}
\includegraphics{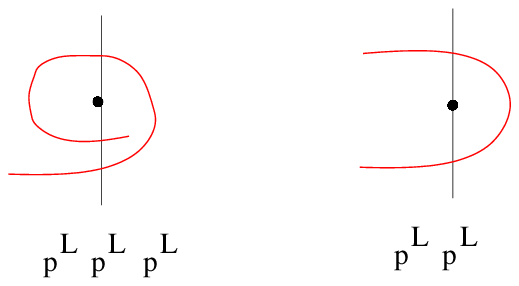}
\end{center}

\noindent 
Since the arc has no 
self-intersections it should wrap around the puncture and then 
unwrap in the opposite direction. In particular it cannot be 
tight. 
\end{proof}

\begin{remark}
One may find however duplicates 
$p^{\epsilon}p^{\epsilon}$, as it can be see in the picture 
above.  
\end{remark}

\paragraph{Untangling braid  terms.} It is known (see \cite{FK2}) that 
$T^*$ is an extension of $T$ by the braid group on infinitely many strands, 
i.e. we have the exact sequence: 
\[ 1\to B_{\infty}\to T^*\to T\to 1 \]
Here $B_{\infty}$ is the group of braids on finitely many 
punctures of $D^*$, i.e. the  ascending union of all braid groups 
$\cup_{n=1}^{\infty}B(D^*_n)$  of finite support. 

\vspace{0.2cm}
\noindent 
There is a natural system of generators for $B_{\infty}$ which was originally 
considered by Sergiescu (\cite{Se}) and then studied by Birman, Ko and Lee 
(\cite{BKL}). For any two punctures associated to adjacent edges 
$e$ and $f$ we associate the braid $\sigma_{ef}\in B_{ \infty}$ which braids 
counter-clockwisely the two punctures and interchange them. Actually 
it suffices to consider only those pairs $(e,f)$ associated to a maximal 
tree in the graph of adjacency of punctures.

\vspace{0.2cm}
\noindent 
Consider now the reduced sequence  $w_{r}(\gamma)$ 
associated to $\gamma$, namely the sequence  
$p_i(\gamma), p_{i+1}(\gamma),...,p_j(\gamma)$ obtained from 
 $w(\gamma)$  by omitting the duplicates i.e. 
we delete $p_{k+1}(\gamma)$ from the sequence 
if $p_k(\gamma)=p_{k+1}(\gamma)$. 
The  consecutive elements in the reduced  
sequence  $w_r(\gamma)$ are still adjacent punctures, lying 
on adjacent edges of the triangulation.  
Suppose that $p_k(\gamma)$ lies on the edge $e_k$. 
We define the {\em untangling braid}, or the {\em  untangling braid  factor}  
$C_{pq}$ by means of the formula: 
\[  C_{pq}=\sigma_{e_ie_{i+1}}^{-\epsilon}\sigma_{e_{i+1}e_{i+2}}^{-
\epsilon}\cdots 
\sigma_{e_{j-1}e_j}^{-\epsilon}\]
\noindent where we put 
\[\epsilon=\left\{\begin{array}{rl}
1, & ~{\rm  if} ~ \epsilon_i=R \\
-1, & ~{\rm  if} ~  \epsilon_i=L
\end{array}\right.\]

\subsection{The existence of conjugate punctures along admissible nonstraight arcs}

\vspace{0.2cm}
\noindent 
Let now $\gamma$ be a tight arc which belongs to the  punctured 
triangulation $\tau^*$ but not to $\tau_B^*$. An arc $\gamma$ with 
the property that there exists a punctured triangulation containing 
it will be called {\em admissible}. As we shall see below, not all 
arcs are admissible. Our algorithm will work only for admissible 
arcs. Furthermore there exists a finite polygon $P\subset \tau_B^*$
which contains all edges from $\tau^*-\tau_B^*$.  
 
\vspace{0.2cm}
\noindent 
Assume that $\gamma$ is uncombed.  We wish to apply the Mosher algorithm 
in order to simplify the arc by means of flips. 
We locate the prong and vertex $v$  from which 
$\gamma$ emerged and denote by $T$ the triangle determined by that 
prong. There are two possibilities: either $\gamma$ is 
combinatorially straight with respect to $T$ or not. 
\begin{enumerate}
\item 
If $\gamma|_T$ is straight then $\gamma$ intersects the edge 
$f$ opposite to 
the vertex $v$. Set $T^*$ for the other triangle of the 
triangulations sharing the edge $f$ with $T$ and  denote 
by $Q$ the quadrilateral $T\cup T^*$. The arc $vv^*$ splits $Q$ into two 
triangles, say $T^L$ and $T^R$ (the superscripts with their 
obvious meaning). 
\begin{enumerate}
\item If $\gamma|Q$ is straight then we use the flip on the edge $f$, 
as in Mosher's algorithm. Thus $f$ is transformed into the edge 
$vv^*$ (with some orientation) and in the new triangulation 
$\gamma$ intersects precisely one triangle $T_2$ 
among $T^L$ and $T^R$. By hypothesis $\gamma|_{T_{2}}$ is again 
straight, and thus we can continue the procedure, as developed below, 
with $T_2$ playing now the role of $T$. 
\item If $\gamma|Q$ is not straight, then let $T_2\in \{T^L,T^R\}$ 
be the triangle containing the prong to which $\gamma$ belongs. 
Then $\gamma|_{T_2}$  is not straight. Notice that $\gamma$ might 
exit $T_2$ and enter next the other triangle. 
However we are now in position to apply the algorithm for the case 
when $\gamma$ is not straight with respect to its first 
triangle it meets. 
\end{enumerate}
\item The case when $\gamma|_T$ is not straight is more involved and it 
will be developed below. 
\end{enumerate}
\noindent We will consider below a way for untangling
arcs which eventually straighten arcs. 
The procedure is based on the following key proposition: 
\begin{proposition}\label{conjugate}
There  is a puncture $q=p_k(\gamma)$ (among those encountered by $\gamma$) 
which is conjugate to $p=p_1(\gamma)$ along $\gamma$. 
\end{proposition}
\begin{remark}
It might happen that $p=q$, but in this case $k\geq 3$. 
\end{remark}
\begin{proof}
Suppose that $\gamma$ emerges from the vertex 
$v$ of the prong with edges $f^L$ and $f^R$. 
By symmetry we can consider that 
$p_1(\gamma)=p(f^L)$. Thus $\gamma$ exits the prong, it crosses 
the edge $f^L$ and goes on the upper halfplane determined 
by $f^L$. Here the halfplane containing the prong was called the 
lower halfplane  and the complementary halfplane the upper halfplane. 

\vspace{0.2cm}
\noindent The proof of this main technical result is given in the next two 
subsections and consists of a detailed analysis of all cases involved. 

\vspace{0.2cm}
\noindent
We say that the arc $\gamma$ is {\em monotone}  it has no conjugate 
punctures. Moreover the arc is {\em L-monotone} if it leaves all punctures 
that it encounters on its left side, except for the puncture 
that it contains. 
 
\vspace{0.2cm}
\noindent Notice that each edge has its endpoints at infinity, or 
alternatively, on the circle at infinity. In particular any edge will 
separate the plane into two  halfplanes. If the edge $e$ is given an 
orientation then it makes sense to consider the halfplanes 
$H_e^L, H_e^R$ which are respectively on the left (or right) of the 
edge $e$.

\subsubsection{The first intersection 
point between $\gamma$ and $f^L$ is different from $p(f^L)$}

\noindent{\bf \em I. The arc $\gamma$ returns on the lower halfplane}

\vspace{0.2cm}
\noindent Suppose next that $\gamma$ crosses again $f^L$ in order to 
arrive in the lower half-plane, leaving all 
punctures that encounters on its left side. Let $\gamma_0$ 
denote the subarc of $\gamma$ travelling in the upper half-plane.  
We will show that these assumptions will lead us to a 
contradiction. 
\paragraph{The puncture $p(\gamma)$ is different from $p(f^L)$.}
Thus $\gamma$ intersects  the edge $f^L$ in a point sitting at 
the right of the puncture $p_1(\gamma)=p(f^L)$, travels around the 
upper halfplane and returns back intersecting again the edge 
$f^L$ in a point still different from  the puncture $p_1(\gamma)$.  
Let the edge $f^L$ have the vertices $v$ and $w$.

\vspace{0.2cm}
\noindent 
The arc $\gamma_0$ has the endpoints on the edge $f^L$ and 
after connecting them by a  line segment we obtain a
circle in the plane. This circle bounds a disk that we call $D$.  
We will use the following well-known lemma going back 
to the proof of Jordan's planar domain theorem: 

\begin{lemma}\label{jordan}
Let $x$ be a point in the plane, which does not belong to 
the boundary $\partial D$ and $[xa$ be a half-line issued from $x$ 
which is transverse to $\partial D$. To each intersection point 
$y\in [xa \cap \partial D$, we associate the number 
$\epsilon(y)\in\{-1,1\}$ which represents the  local 
algebraic intersection number between $[xa$ 
(whose orientation points towards $a$)
and  $\partial D$. Then 
\[ \left|\sum_{y\in [xa \cap \partial D} \epsilon(y)\right| \in\{0,1\}\] 
Moreover, this value is independent on the choice of the 
half-line $[xa$. 
Furthermore, the claim holds true when  $x\in \partial D$. 
\end{lemma}
\begin{proof}
Actually the sum is $0$ if $x\not\in D$ and $1$, if $x\in D$.
When $x\in\partial D$ we deform slightly $x$ off $a$. Then 
the value for $x$ is the same as that associated to the 
perturbed point.   
\end{proof}

\begin{lemma}
Consider now a triangle $T$ made of ideal arcs which has 
non-empty intersection with the curve $\gamma_0$. Then the 
homeomorphism type of the pair $(T, T\cap \gamma_0)$ 
belongs to one of the following patterns 

\begin{center}
\includegraphics{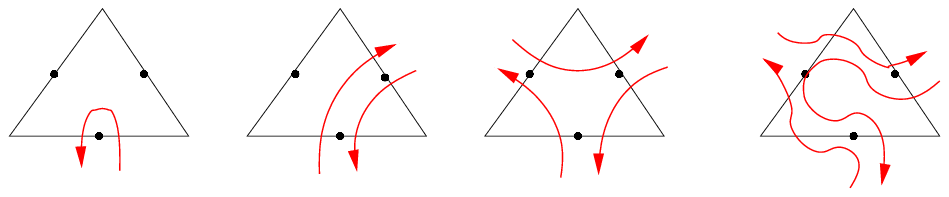}
\end{center}
\end{lemma}
\begin{proof}
Let us see the consequences of  
lemma \ref{jordan} when applied to a half-arc determined by the puncture 
$p(e)$ of some edge $e$ which is crossed by $\gamma_0$. 
Assume that the half-arc is crossed at least once by $\gamma_0$. 
Assume furthermore that $p_e$ does not belong to $\gamma_0$. 
According to our hypothesis all crossing points should leave the 
puncture $p_e$ on their left. This implies that the local
algebraic intersection number at a crossing point between $\gamma_0$ and 
the half-arc (oriented towards 
the vertex at infinity) has always the same value, and in our particular 
situation where $p_e$ is on its left side, it  should be positive.  
Now, the lemma above implies that 
we must have only one crossing point between the half-arc and 
$\gamma_0$, since otherwise their sum would be greater than $1$.  
Moreover, if there is a crossing point between one 
half-arc of the edge $e$ and $\gamma_0$ then we should have at least 
one intersection point between $\gamma_0$ and the other half-arc 
of $e$ issued from the puncture $p(e)$. This follows from the second 
part of the lemma above. In particular, we obtained that any edge 
is crossed twice, each half-arc being crossed precisely once.

\vspace{0.2cm}
\noindent
Another possibility to take into account is when $p_e$ belongs to 
$\gamma_0$. Assume that the arc $e$ is not tangent to $\gamma_0$. 
If some half-arc determined by $p_e$ is crossed at least 
twice  by $\gamma_0$ then deform slightly $p_e$ 
towards the vertex at infinity. 
We will obtain a point for which the half-line which it determines 
has algebraic intersection at least $2$, by the same argument as 
above. This contradicts  the lemma. Thus each half-arc can be crossed 
by $\gamma_0$  at most once more. We claim that only one  half-arc 
among them can have nontrivial intersection with $\gamma_0$. 
Assume the contrary, namely that both half-arcs $p_ea$ and $p_eb$ 
intersect $\gamma_0$. 
We have that the algebraic intersection number of 
$\gamma_0$ and one half-arc (oriented towards infinity) 
is positive. Thus
the local algebraic intersection number of $\gamma_0$ and $p_ea$ 
should be negative at $p_e$, otherwise their sum being at least 2. 
Thus the frame $(p_ea,\dot\gamma_{p_e})$ 
is negatively oriented. Similarly, the frame 
$(p_eb, \dot{\gamma}_{p_e})$ should be negatively oriented, which 
is impossible because these two frames have opposite orientations. 
This shows that the arc $\gamma$ intersects each edge precisely twice, 
with the possibility that one intersection point be the puncture. 

\vspace{0.2cm}
\noindent There is one more possibility, when the arc $\gamma_0$ is 
tangent at the puncture $p_e$ to the edge $e$. The argument above 
shows that  in addition to what we already 
saw $\gamma_0$ can intersect once more each half-arc 
$p_ea$ and $p_eb$ when $p_e$ is a tangency point. The local model 
is that from below: 

\begin{center}
\includegraphics{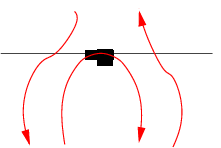}
\end{center}

\vspace{0.2cm}
\noindent Consider now a triangle  $T$ that $\gamma_0$ intersects. 
There are only finitely many possibilities for $\gamma_0\cap T$, 
so that any edge is crossed precisely twice and 
there are no self-intersections, and these are precisely those pictured 
above. 
\end{proof}


\noindent We are able now to formulate the following result which 
explains the form of the curve $\gamma_0$. 
Let $\Gamma$ be a finite  subtree of the 
tree of the triangulation $\tau_B$, having the puncture 
$p_1$ as one of its leaves. If there is a tangency puncture between 
$\gamma_0$ and an edge then consider that in the dual graph we add 
edges between the punctures of the respective triangle and remove 
the associated $Y$ graph. 
Construct the  planar regular neighborhood $N(\Gamma)$ of $\Gamma$ and
 consider its boundary  
$\partial N(\Gamma)$.

\begin{lemma}[Disk lemma]\label{tree}
The curve $\gamma_0$ is the intersection of 
$\partial N(\Gamma)$  with the upper-half plane determined by $f^L$, 
which is $N(\gamma)$ minus a small cap around 
$p$.  
\end{lemma}
\begin{proof}
Remark that the arc $\gamma_0$  might pass thru a puncture. 
Each model above (and its images under the $\Z/3\Z$ symmetries) 
can appear within $\gamma_0$. Moreover, we can now obtain $\gamma_0$ 
using this Lego toolkit by gluing up triangles with the models 
inside, which have matching boundaries. The arcs that we obtain 
are described as in the statement. 
\end{proof}

\vspace{0.2cm}
\noindent We are able now to finish the proof of the proposition 
in the 
case under scrutiny. In fact, we obtained that the arc $\gamma_0$ 
(hence $\gamma$) comes back crossing again $f^L$ in one point 
which belongs to the half-edge $vp$.
If the crossing point is not the puncture, then  $\gamma$ should enter 
the domain determined by $vp$ and the first part of arc sitting in 
the lower half-plane. This domain does not contains any other vertex 
of the triangulation, and thus $\gamma$ has to exit the domain in 
order to abut to some vertex (different from $v$). 
But $\gamma$ cannot 
cross itself and $\gamma$ cannot cross the edge $vw$ again since 
all edges are crossed twice. This is a contradiction. 

\begin{center}
\includegraphics{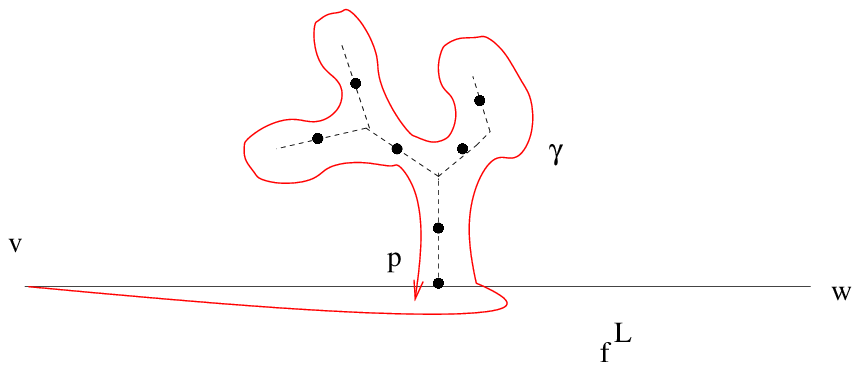}
\end{center}

\paragraph{The arc $\gamma$ return on $f^L$ and hits it at $p_1$.}
Another possibility left is that the second intersection point 
between the arc $\gamma$ and the edge $f^L$ coincides with the 
the puncture $p$. 
\begin{itemize}
\item If the arc $\gamma$ enters the lower half-plane then 
we get a contradiction by the same argument that we used above. 
\item Otherwise, the arc $\gamma_0$ is tangent at $f^L$ at $p$ and 
goes up in the upper half-plane. There are  again 
two possibilities:
\begin{itemize}
\item The orientations of the two tangent arcs are compatible. 
Then the arc $\gamma_0$ should cross once more the segment $pw$ and 
thus the puncture $p$ will be on its right side.  
\item The two tangent arcs point in opposite directions, which we suppose 
to be the case from now on. 
\end{itemize}
\end{itemize}
We denote by $\gamma_1$ the subarc of $\gamma$ issued from the 
puncture $p$ and lying in the upper half-plane.  
Let  $f^L$ be the edge of the triangle 
$uvw$ sitting in the upper half-plane. We have several  possibilities, as 
could be seen from the picture below: 

\begin{center}
\includegraphics{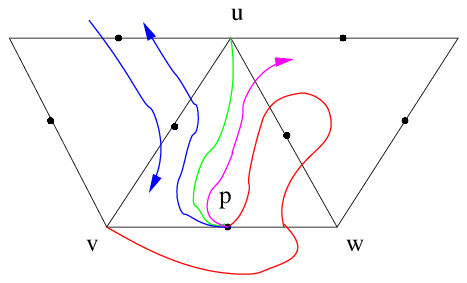}
\end{center}

\begin{enumerate}
\item If the arc $\gamma_1$ crosses first the edge $uw$. 
Since $\gamma_1$ is disjoint from $\gamma_0$ it should 
cross $vw$ leaving the puncture on its right side. 
\item If $\gamma_1$ goes straight to the vertex $u$. 
We shall see below that in this case the arc $\gamma$ cannot be 
admissible. 
\item Otherwise $\gamma_1$ crosses first the edge $uv$.
\begin{enumerate}
\item If $\gamma_1$ crosses again the edge $uv$, then 
the same argument which was used for the edge $f^L$ applies. 
Thus $\gamma_1$ returns on $uv$ by intersecting it 
in a point of the  half-arc containing $v$. This time 
$\gamma_1$ should avoid the puncture. Further, $\gamma_1$ is now 
forced to reach the vertex $v$ (which is a contradiction) or else 
to enter the region $N(\Gamma)$ containing no vertex, either. 
Once more the arguments above show that 
$\gamma_1$ should turn back and cross $pw$ and thus $p$ is on its right side.  
\item Eventually, if $\gamma_1$ reaches a vertex in the left upper
half-plane determined by $uv$ then we will see that $\gamma$ cannot be 
admissible. 
\end{enumerate} 
\end{enumerate}
\noindent In order to deal with these cases we will open a parenthesis 
in the proof in order to state an intermediary result, before resuming.

\paragraph{The admissibility lemma.}
The technical result below will be used several times in the sequel. 

\begin{lemma}[Admissibility lemma]\label{admissible}
Let $\alpha$ and $\beta$ be two admissible oriented arcs (coming from 
possibly different triangulations) with the same endpoints. 
Let denote by $X(\alpha,\beta)$ the set of those punctures 
which are on the left hand side of $\alpha$ but on the 
right hand side of $\beta$. Similarly for $X(\beta,\alpha)$. 
Then the number of elements of $X(\alpha,\beta)$ is the same 
as the number of elements of $X(\beta,\alpha)$. 
\end{lemma} 
\begin{proof}
Since $\alpha$, $\beta$ are admissible there exist 
punctured triangulations 
$\tau_1,\tau_2$ containing them. Moreover, these triangulations are 
identical outside some polygon $P$ and they coincide with $\tau_B^*$. 
Thus there exists an element 
$\zeta$ of $T^*$ which sends $\tau_1^*$ with the d.o.e. $\alpha$ 
onto $\tau_2^*$ with the d.o.e. $\beta$. This means that 
there exist  triangulated polygons $P_j$ which 
are subpolygons of $\tau_j$ such that $\tau_1-P_1$ 
is rigid homeomorphic to $\tau_2-P_2$. We can choose 
$P_j$  large enough in order to contain $P$ in its interior.
The arc $\alpha$ (respectively $\beta$)  splits $P_1$ 
(respectively $P_2$) into its left part $P_1^L$ (and $P_2^L$) and 
its right side $P_1^R$ (and $P_1^R)$. 
Let us order circularly the edges of $P_1^L$ as 
$a_1,a_2,...,a_k, \alpha$ and those of $P_2^L$ as 
$b_1,b_2,...,b_m, \beta$. 
Since $\tau_1-P_1$ and $\tau_2-P_2$ are rigid equivalent 
and the arcs $\alpha$ and $\beta$ correspond to each other, it 
follows that the edges $a_i$ and $b_i$ should correspond to each 
other by means of this rigid homeomorphism. 
Further, there are no vertices at infinity at the interior 
of $P_j$ and thus the only possibility to arrange the 
convex polygons $P_i$ in the plane is as in the picture below
(where $P_j^R$ are drawn in dotted lines), namely:  
$a_1$ is surrounded by the edges $b_1,b_2,...,b_{i_1}$
$b_{i_1+1}$ is surrounded by the edges $a_2,a_3,...,a_{j_1}$ etc. 
In particular, there exists a polygon $P_1\cap P_2$ which is admissible 
(and thus embedded into $\tau_B^*$) which has the following 
edges in the left hand side of  
the arcs $\alpha,\beta$: $a_1,b_{i_1+1},a_{j_1+1},...,b_k$.
In particular one finds as many edges from $P_1$ as edges  from $P_2$, 
and it might happen that some edges of $P_1\cap P_2$ 
are common to both if some $a_i$ equals some $b_j$.

\begin{center}
\includegraphics{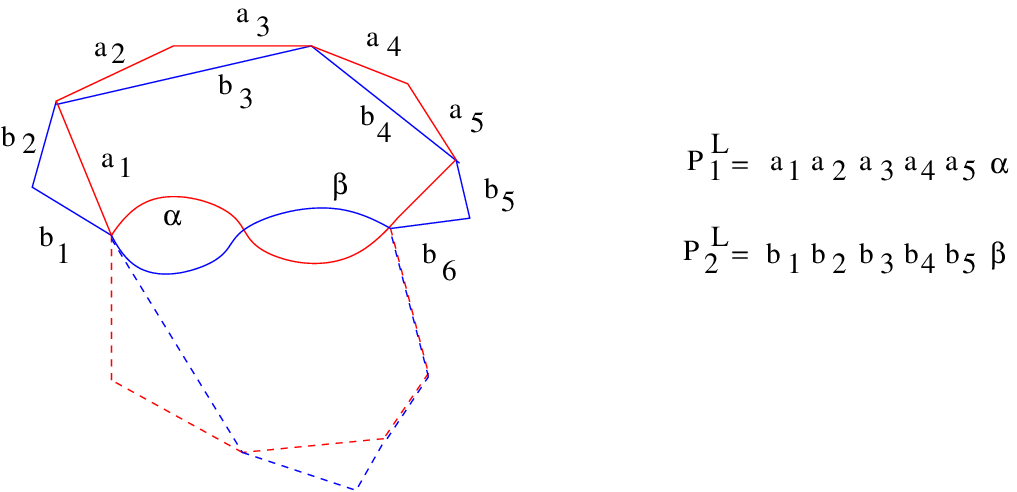} 
\end{center}

\noindent The polygon $Q=P_1\cap P_2$ is split by $\alpha$ into 
Let $Q_{\alpha}$ be the part of the polygon $P_1\cap P_2$ 
sitting on the left of $\alpha$ and similarly $Q_{\beta}$. 
Thus $Q_{\alpha}$ and $Q_{\beta}$ have in common all edges 
but the arcs $\alpha,\beta$. The common part of their 
boundaries is made of $r$ edges from $P_1$ and $r$ edges from $P_2$ 
(some edges might possibly belong to both).

\vspace{0.2cm}
\noindent 
Recall further that the element $\zeta\in T^*$ should send 
homeomorphically $P_1^L$ onto $P_2^L$, by sending 
boundary edges to boundary edges and the arc $\alpha$ onto $\beta$. 
In particular, the number of punctures into $P_1^L$ should 
be equal to the number of punctures in $P_2^L$. 
Moreover, each polygon with $n$-vertices  should contain 
$n-2$ boundary punctures
and $n-3$  interior punctures. Since $Q_{\alpha}$ and $Q_{\beta}$
contain as many edges from $P_1^L$ as edges from $P_2^L$ 
we find that the  union of all polygons from 
$P_1^L-Q_{\alpha}$ contains the same number of 
punctures as the union of all polygons from 
$P_2^L-Q_{\beta}$. In fact all these polygons are admissible 
subpolygons of $\tau_B^*$ and their total number of edges is 
the same in both cases. As a consequence, the number 
of punctures within $Q_{\alpha}$ coincide with the number 
of punctures in $Q_{\beta}$. This is implies the 
statement of the lemma. 
\end{proof}

\begin{remark}
It is useful to be more precise 
concerning the points which are on the left 
of $\alpha$ and the right hand side of $\beta$. 
The (oriented) arc $\gamma$ has endpoints which are 
vertices at infinity, 
meaning that they lay on the boundary circle of the compactification 
disk of the hyperbolic plane. These endpoints determine an arc 
of circle lying at the left of $\gamma$. The union of this 
arc  at infinity with $\gamma$ is a circle bounding the disk 
$D^L(\gamma)$. The complementary disk is $D^R(\gamma)$ sitting on 
the right of $\gamma$.  A point is said to be at the left 
of $\alpha$ and  at the right of $\beta$ if it belongs to 
$D^L(\alpha)-D^R(\beta)$.  
\end{remark}

\paragraph{End of the proof when the arc $\gamma$ return on $f^L$ and hits it at $p_1$.} We are able now to finish the proof in the cases 
enumerated above. Namely, the arc $\gamma$ is made of two pieces, 
one being $\gamma_0$ which crosses $f^L$ and goes in the upper plane 
then hits $f^L$ at the puncture $p$; then the second piece 
$\gamma_1$ emerges into 
the  left upper half-plane determined by $uv$ abutting to some 
vertex  $z$ there (which might be $u$). Our claim is that 

\begin{lemma}
If $\gamma_1$ is L-monotone then $\gamma$ cannot be admissible. 
\end{lemma}
\begin{proof}
The main tool is the admissibility lemma \ref{admissible}. 
The first step is to construct an admissible arc $\delta$ joining 
$v$ to $z$. Consider the dual graph associated to the 
triangulation and next the geodesic joining the
(vertex dual to the) triangle opposite to $uvw$ and having the 
edge $uv$ in common with it, to the closest triangle having 
$z$ among its vertices. Passing again to the 
dual, the union of triangles corresponding to 
vertices of this geodesic is the polygon $M(v,z)$. 
This is the smallest polygon made of adjacent triangles joining 
$v$ and $z$. We can realize $M(v,z)$ as a convex polygon in 
the plane which is triangulated by means of several 
diagonal edges.  Moreover the line segment
$vz$ has to intersect all diagonal edges, since one could get 
rid of any triangle disjoint from this line segment, which would 
contradict the minimality. 

\vspace{0.2cm}
\noindent One can use iteratively flips within the polygon 
$M(v,z)$ in order to weakly comb $zv$, which means that 
we find an arc joining $v$ and $z$ which remains within the 
polygon. There are several possibilities, by making inductively 
the diagonal edges to turn from $v$ closer and closer to 
the vertex $z$. We are constrained to deal only 
with arcs passing thru the puncture of the quadrilateral. Thus 
the final arc that we obtain is neither canonical nor the line 
segment, but it will be convenient for our purposes. We denote 
any such arc by $\widetilde{vz}$. 

\vspace{0.2cm}
\noindent 
We will show that there cannot be any puncture which sits 
at the left of $\widetilde{vz}$ and at the right of the arc 
$\gamma_1 \cap M(v,z)$. 
The arc $\gamma_1$ enters $M(v,z)$ throughout the right half-arc of 
the edge $vu$. Moreover, $\gamma_1$ should cross all 
diagonal edges of $M(v,z)$ since their endpoints are on the circle 
at infinity and they separate $v$ from $z$. Also $\gamma_1$ cannot 
hit any other puncture since $\gamma$ already contains $p$. 

\begin{lemma}
The arc $\gamma_1$ hits any diagonal edge $e$ precisely once. 
\end{lemma}
\begin{proof}
Assume that $\gamma_1$ hits $e$ once at $q$ and 
then returns and intersects 
$e$ again. We claim that any further intersection point of 
$\gamma_1$ with some edge $f$ already crossed by $\gamma_1$ should be 
closer to the puncture $p_f$ than the previous hit. In particular 
$\gamma_1$ should cross again $e$ (since $e$ separates $v$ from $z$)
and it will cross it by entering throughout  the segment $p_eq$. 
Now,  the subarc of $\gamma_1$ between the first and the second hit 
of $e$ union with a small segment on $e$ will  bound a disk
which is disjoint from $z$. Thus the arc $\gamma_1$ cannot escape 
this disk (by the Disk lemma \ref{tree}) and thus cannot reach $z$, 
which is false. 
  
\begin{center}
\includegraphics{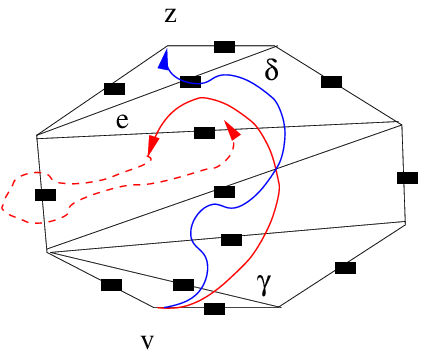}
\end{center}

\end{proof}
\noindent Moreover, if $\gamma_1$ hits any diagonal edge $e$ 
precisely once then the puncture $p_e$ belongs to $D^L(\gamma_1)$. 
This happens thus for all punctures from the interior of 
$M(v,z)$.  Now, we will consider the punctures from 
$D^L(\delta)$ which are not interior points of $M(v,z)$. 
Such a point cannot be in $D^R(\gamma_1)$ unless 
$\gamma_1$ crosses one boundary edge $f$ of $M(v,z)$ 
located at the left of $\delta$. Since $z$ is in $M(v,z)$, 
the arc $\gamma_1$ should cross again $f$. The arguments above 
show that $\gamma_1$ crosses once more $f$ on the other half-arc. 
The lemma above implies that, in order to reach $z$ 
the arc $\gamma_1$ should either cross itself or else 
cross again $f$. But when crossing $f$ again the arc will enter 
again a disk and  by the Disk Lemma \ref{tree} it cannot escape anymore.  This proves that 
$\gamma_1$ cannot cross the left part of the boundary of 
$M(v,z)$ and thus the set of punctures in the left of $\delta$ and 
the right of $\gamma_1$ is empty. 
\end{proof}

\vspace{0.2cm}
\noindent We return now to  the arc $\gamma$ which contains the
extra  piece $\gamma_0$. We saw that $\gamma_0$ contains at its left 
at least one more puncture out of $p$. Moreover, all punctures 
from the lower half plane determined by $uv$ are from $D^R(\delta)$.   
This implies that $X(\gamma,\delta)$ has at least one element, while 
$X(\delta,\gamma)$ is empty. This contradicts the admissibility lemma. 
This proves therefore the claim in the case under consideration. 

\vspace{0.2cm}
\noindent

\noindent{\bf \em II.The arc $\gamma$ remains in the upper half-plane}

\vspace{0.2cm}
\noindent
From now on we will consider the situation when  
$\gamma$ remains in the upper half plane and therefore reaches a vertex $z$ 
in the upper half-plane, which might coincide 
with the  other vertex $w$ of the edge $f^L$. 
There are again two possibilities for $z$: either $z$ lies within the 
right upper half-plane determined by $uw$ or else in the left 
upper half-plane determined by $uv$, the point $u$ inclusively. 
The second case is very similar to the last part above. 

\vspace{0.2cm}
\noindent
Assume that $z$ belongs to the 
right upper half-plane determined by $uw$. 
Then either $\gamma$ crosses first  $uw$ or else $uv$, or it goes 
towards $u$. 

\vspace{0.2cm}
\noindent
If $\gamma$ crosses $uv$ then it has to return back and cross 
$uv$ again in order to arrive at $z$. The Disk lemma \ref{tree}  shows that 
$\gamma$ cannot escape from the left upper half-plane. 

\vspace{0.2cm}
\noindent
If $\gamma$ goes to $u$ then the admissibility lemma \ref{admissible}
shows that the arc is not admissible. 

\vspace{0.2cm}
\noindent
Thus $\gamma$ crosses $uw$. If $\gamma$ crosses again $uw$ 
then the disk argument shows that it cannot escape towards a vertex. 
Hence $\gamma$ has only one intersection point with $uw$. 

\vspace{0.2cm}
\noindent
Let now $M(v,z)$ be the polygon with the smallest number of triangles 
which join $v$ to $z$. We already saw above that it should exist 
an admissible arc $\delta$ joining $v$ to $z$ and lying within 
$M(v,z)$. We can be more precise, as follows:

\begin{lemma}
There exists an admissible arc $\delta$ joining 
$v$ to $z$ within $M(v,z)$ such that the number of punctures 
encountered by $\delta$ which are on its right side is equal to 
$N-1$, where $N$ is  computed as follows. Let us consider the two 
vertices of the dual graph associated to the triangle $uvw$ and 
the triangle $T_z$ containing $z$ within $M(v,z)$. Join the two 
vertices by a geodesic in the binary tree. Then $N$ denotes the 
number of times the geodesic has to turn left at some intermediary vertex. 
Here we assume that the first segment of the geodesic issued 
at $uvw$ turns right, because the next triangle is adjacent to $uw$.   
\end{lemma}
\begin{proof}
We will change the triangulation of the polygon $M(v,z)$ in order 
to connect $v$ to $z$. Remark that $M(v,z)$ can be seen as a convex polygon 
in the plane. Let denote the vertices of $M(v,z)$ in clockwise order 
$v,u=u_1,u_2,...,u_k,z,w_m,w_{m-1},...,w_1=w$. 
The vertices $u_j$ will be called upper vertices and 
the vertices $w_j$ lower vertices. 
We  will change inductively the diagonals  edge by means of flips 
aiming at combing the segment $vz$. 
This means that we will adjoin step by step all intermediary diagonals 
$vu_2,...,vu_k$ and intercalated among them 
$vw_2,...,vw_m$. The new diagonals obtained by flipping will be not touched 
on the next steps and will be called new arcs. At the end we will get the 
triangulation containing all $vu_j, vw_j$ and $vz$. 
There is only one issue to care about: when the diagonal arc  
$vu_j$ has been adjoined by means of a flip into some former 
diagonal edge $u_jw_k$, the new arc $vu_j$ was constrained 
to pass thru the puncture of the diagonal  $u_jw_k$ 
crossed at the previous step. 
Thus the location of the puncture associated to such a diagonal arc 
is determined by the time when the new diagonal arc is adjoined. 

\vspace{0.2cm}
\noindent
Now, a flip on $u_1w_1$ will transform it into either 
$vu_2$ or else into $vw_2$, depending on whether 
$w_1$ was adjacent to $u_2$ or else $u_1$ to $w_2$. 
If the former case happen then 
the puncture $p_{u_1w_1}$ will be said to be an upper puncture,
otherwise it will be called a lower puncture. 
In general, upper punctures will belong to upper diagonals 
$vu_j$ and lower punctures to lower diagonals $vw_j$. 
Assume that at some step we get an 
upper puncture on some arc  $vu_j$ after flipping the 
former diagonal edge $u_jw_k$. All new arcs which will be adjoined  
at the next steps, namely $vu_i$ with $i\geq j+1$ and 
$vw_l$ with $l\geq k+1$ will cross the former edge $u_jw_k$, and they 
should leave that upper puncture on their left side. 
This is so because otherwise they would  
cross the arc $vu_j$ just adjoined. 
Similarly, lower punctures will remains always to the right 
of the new arcs. 

\begin{center}
\includegraphics{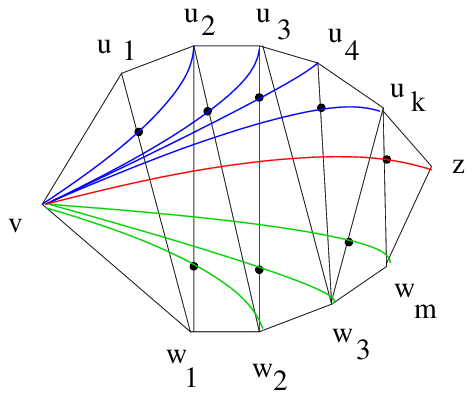}
\end{center}

\noindent 
At the end of our combing we will have then a number of 
$k$ upper punctures associated to 
the upper diagonal arcs $vu_j$ and also a number of $m$ lower punctures 
associated to the lower diagonal arcs $vw_j$. Thus the number of punctures 
encountered by $vz$ which rest on its right side is $m$. 
Eventually the value of $m$ can be easily checked to be that claimed 
by recurrence on the number of triangles involved in $M(v,z)$. 
\end{proof}
\begin{corollary}
If $d(z,w)=r\geq 1$ then there exists an admissible arc $\delta$ joining 
$v$ to $z$ within $M(v,z)$ such that the number of punctures 
encountered by $\delta$ which are on its right side is at least $r-1$.  
\end{corollary}
\begin{proof}
This follows from the fact that $N\geq d(z,w)$, because there exists 
a path $w_1w_2...w_Nz$ joining $w$ to $z$ of length $N$. 
Notice that in general the distance $d(z,w)$ might be smaller than 
$N$ computed above, because there might exist a shorter path 
using upper vertices.  
\end{proof}

\vspace{0.2cm}
\noindent Now the arguments are similar to those from the case when 
$\gamma$ returns to the lower halfplane. 
If $\gamma_0\subset M(v,z)$ then observe that any diagonal edge should 
be crossed precisely once, otherwise we can find a disk containing the 
arc and  by the Disk  lemma \ref{tree}  
the arc would not escape from it. This implies that 
there is no puncture in the interior of $M(v,z)$ which lies in 
$D^R(\gamma)$. There exists only one puncture which is on $\gamma$. 
Further, the puncture $p=p_{vw}$ should belong to $D^L(\gamma)$. 
The admissibility lemma implies that the arc $\gamma$ is not admissible 
as soon as $N\geq 1$, which is implied by $d(z,w)\geq 1$.  
The same argument shows this is the case also when $\gamma_0$ crosses 
boundary edges of $M(v,z)$ on the right. On the other hand 
if $\gamma$ crosses a boundary edge at the left of $\delta$ then 
the Disk lemma \ref{tree}  will lead us to a contradiction. 

\vspace{0.2cm}
\noindent There is one more case left, namely when $z=w$. In this case 
the arguments above break down. However, we can take $\delta=vw$. 
Further, it should not exist any other puncture encountered by $\gamma$ on 
its right side and thus $\gamma$ is isotopic to the following configuration: 

\begin{center}
\includegraphics{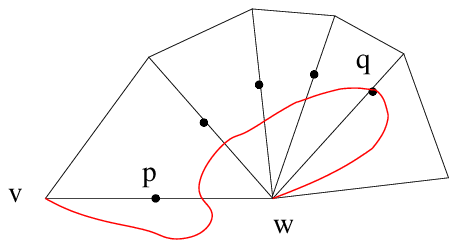}
\end{center}

\vspace{0.2cm}
\noindent This means by definition that the puncture $q$ (lying on $\gamma$)
is conjugate to $p$.

\subsubsection{The first 
intersection point between $\gamma$ and $f^L$ is  $p(f^L)$}
\vspace{0.2cm}
\noindent
Assume that there are no conjugate points along $\gamma$. 
We also suppose that $\gamma$ is tight and thus there are no 
superfluous tangencies. Then 
$\gamma$ should be L-monotone. In fact, any point which lives on the 
right side of $\gamma$ is conjugate to $p=p_1(\gamma)$, by definition.
Further the arc $\gamma^{-1}$ which is $\gamma$ with the reverse 
orientation should be R-monotone. Assume that there is at least one 
puncture that is encountered by $\gamma$ out of $p$. Then according to 
the Part I of the proof of proposition \ref{conjugate} that puncture 
should be conjugate to another puncture along $\gamma^{-1}$ and this 
would contradict the fact that $\gamma^{-1}$ is monotone. 
Thus $\gamma$ cannot meet any other puncture and thus the arc is 
isotopic either to $f^L$ or else to the arc below: 

\begin{center}
\includegraphics{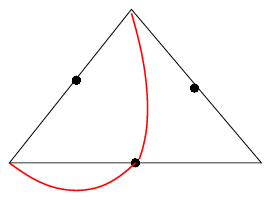}
\end{center}
 
\vspace{0.2cm}
\noindent However it can be easily seen that this arc is not admissible 
by means of the admissibility lemma. This establishes the claim. 
\end{proof}

\subsection{Simplifying arcs by combing and straightening tight arcs}

\subsubsection{Combing admissible arcs}
\noindent Denote  by 
$L(\gamma)$ the length of $\gamma$, namely the 
total number of punctures encountered by $\gamma$. 

\begin{lemma} 
Recall that $p_1(\gamma)$ and $p_k(\gamma)$ are conjugate 
along $\gamma$. Then the arc 
$\gamma_1=C_{p_1(\gamma)p_k(\gamma)}\gamma$
has length $L(\gamma_1) < L(\gamma)$. 
\end{lemma}
\begin{proof}
Suppose that $\gamma$ emerges from the vertex 
$v$ of the prong with edges $f^L$ and $f^R$. 
By symmetry we can consider that 
$p_1(\gamma)=p(f^L)$ and thus $\gamma$ crosses 
$f^L$ and goes on the upper plane determined by $f^L$. 
There are a number of punctures encountered by $\gamma$ 
and left on its left side until it reaches $p_k$ on its 
right side. 

\begin{lemma}
The puncture $p_k$ belongs to the upper half plane or 
it coincides with $p_1$. 
\end{lemma}
\begin{proof}
Otherwise $\gamma$ crosses again $f^L$ in order to 
arrive in the lower half-plane, leaving all 
punctures encountered on its left side. 
If the new crossing point is on the right of the 
previous crossing point then $\gamma$ is not tight, 
as already observed in the proof of lemma \ref{conjugate}. 
If the second cross point lies on the segment $vp_1$ 
then the arc $\gamma$ unwraps around $p_1$ but it has 
to unwrap and exit the same way, and thus it cannot be tight. 
The remaining possibilities are that $\gamma$ 
crosses again $f^L$ leaving $p_1$ an the right hand side, or 
else that $p_k$ is on the upper half-plane.  
\end{proof}
\noindent We can now verify that the 
untangling transformation $C_{p_1p_k}$
permutes the punctures
by bringing $p_k$ on the position of $p_1$ and translates 
each other $p_j$ onto the next $p_{j+1}$. Further, the 
image of the arc $\gamma$ by means of the untangling 
braid  has the following shape: 
\begin{center}
\includegraphics{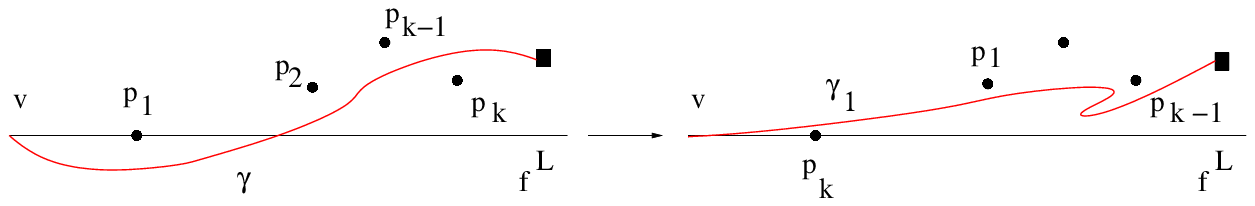}
\end{center}
Thus $\gamma_1$ avoids the prong determined by $f^L$ and $F^R$ 
and belongs to the prong sitting on its left side (having 
$f^L$ as edge). All punctures  $p_k$ but $p_1$ are 
still encountered by $\gamma_1$. 
Moreover, the mapping class $C_{p_1p_k}$ can be 
represented by a homeomorphism whose support is 
contained in a neighborhood of the graph of vertices 
$p_j$ and edges $p_jp_{j+1}$. In particular, the arc $\gamma_1$ 
does not intersect any other edge of the triangulation which 
was not previously crossed by $\gamma$, because each 
$p_j$ corresponds to an edge and the segment $p_jp_{j+1}$ 
lies in the interior of a triangle. This proves that 
$\gamma_1$ encounters precisely the punctures 
located at $p_2,p_2,...,p_k,...$ and avoids $p_1$. 
This proves that $L(\gamma_1)<L(\gamma)$.   
\end{proof}

\noindent 
The arc $\gamma_1$ is then simpler that $\gamma$. 
If $\gamma_1$ is not yet straight with respect to its 
first triangle that it meets then we apply again lemma
\ref{conjugate} and use the associated untangling 
transformation in order to obtain the arc $\gamma_2$. 
We continue this way by defining recurrently 
the arcs $\gamma_{m+1}$ obtained by untangling $\gamma_m$. 
\begin{lemma}
After finitely many steps the arc $\gamma_m$ is straight 
within the first triangle that it meets. 
\end{lemma}
\begin{proof}
At any step the length $L(\gamma_{m+1})< L(\gamma_m)$. 
In particular for $m$ large enough $L(\gamma_m)=1$.
This means that $\gamma_m$ does not encounter any other 
punctures but the puncture that it contains. This puncture 
could belong to the opposite side of the vertex $v$, or 
else on a side of the prong containing $\gamma_m$.
The arc $\gamma$ continues after the puncture and 
reaches a vertex without crossing any other edge.
This implies that the arc goes straight to the opposite vertex 
or else goes along one  half-arc in an edge. 
Summing up, the arc $\gamma$ has one of the shapes 
pictured in the figure below, namely: 
$vpv^*$, $vqt$, $vpz$, $vpt$, $vqz$, $vqw$. 
The first two arcs are straight. 

\begin{center}
\includegraphics{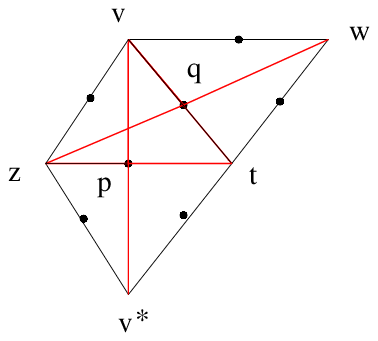}
\end{center}
\noindent We claim that none of the remaining arcs
is admissible. Let us concentrate on the arc $vpz$. 
Since $zpt$ is an edge of the triangulation 
obtained so far $\widetilde{\tau^*}$, it is 
an admissible arc. Recall that the arc $\gamma$ was 
an admissible arc within the  finite polygon $P$ with the 
property that all arcs outside $P$ belong to $\tau_B^*$. 
All flips and  untangling transformations  used up now are 
supported in $P$, which means that they do not touch any arc 
outside $P$. If $vpz$ were an arc of some triangulation, this 
means that it is an arc of the triangulation obtained from 
$\tau^*$ by means of some flips and untangling transforms 
which aimed at combing it. Thus $vpz$ would be an arc 
of a triangulation $\widetilde{\tau}^*$ which differs from 
$\tau_B^*$ only within the polygon $P$. 
It makes sense therefore to consider the element of $T^*$ 
which sends $\widetilde{\tau^*}$ with the d.o.e. 
$zpt$ onto the triangulation $\widetilde{\tau^*}$ with the d.o.e. 
$zpv$.  Recall now that elements of $T^*$ are mapping classes 
of homeomorphisms supported in some finite polygon, which in this 
case is $P$. Further, any ideal arc $\eta$  divides the 
compactified plane into two disks $D(\eta)^+$ 
and $D(\eta)^-$, each disk intersecting 
the polygon $P$ inside a disk. Further a homeomorphism 
of $P$ sending one ideal arc into another one should 
send $D(zpt)^+\cap P$ onto $D(zpv)^*\cap P$ and punctures 
onto punctures. This is a 
contradiction  because the number of punctures of 
$D(zpt)^+\cap P$ is at least one unit greater than the number of 
punctures inside $D(zpv)^+\cap P$ since the former 
disk should contain the puncture $q$. This shows that 
$zpv$ cannot be admissible. The same argument works for the 
other arcs. The only possibility for $\gamma_n$ is 
to be either an edge $vqt$ or a diagonal $vpv^*$, and thus 
straight. 
\end{proof}

\subsubsection{Straightening combed arcs}

\noindent 
Let us consider the case when the arc $\gamma$ is combed but 
it is not isotopic to the corresponding arc $xy$ determined by its 
endpoints $x,y$ in $\tau_B^*$. 

\vspace{0.2cm}
\noindent 
Suppose that $\gamma$ belongs to a prong $Q_1$ contained in the 
triangle $T_1$. 
\begin{enumerate}
\item If $\gamma|{T_1}$ is not straight, then use the procedure 
from the previous section and find an untangling  braid term that makes it straight. 
\item If  $\gamma|{T_1}$ is straight then,
\begin{enumerate}
\item If $T_1$ is one of the two triangles containing $xy$ then 
$\gamma$ is isotopic to $xy$. 
\item Otherwise let $f_1$ be the edge opposite to the prong  $Q_1$ to which 
$\gamma$ belongs. Assume that we do a flip on $f_1$, which changes 
the base triangulation but will {\em not} be recorded in the 
combing word. The effect of the flip is to split the prong $Q_1$ into two 
prongs of the new triangulation. We assume that $Q_2$ is the new prong 
to which $\gamma$ belongs and $T_2$ is the triangle that it contains it. 
Then we iterate the procedure above. 
\end{enumerate}
\item The procedure stops when $\gamma$ has been straightened and 
is isotopic to $xy$. 
\end{enumerate}
\vspace{0.2cm}
\noindent The word  associated to the straightening procedure is 
the product of all untangling terms used. This is an element 
of $B_{\infty}$ since $\gamma$ is already combed.

\subsection{Complements on straightening arcs}
\noindent
There is also a {\em global} straightening procedure  which works 
for  any subarcs  (not only initial ones) of a given arc and for  
combed arcs as well. 
Suppose that  $\gamma$ has the endpoints $\gamma(0)$ and 
$\gamma(1)$ which determine the  line segment 
$[\gamma]=\gamma(0)\gamma(1)$ which we call the shadow of $\gamma$. 
If $\gamma$ is combed  then its shadow is an edge of $\tau_B^*$. 
If we are looking only to a subarc of $\gamma$ which has to be 
straightened  we will compare it with its shadow. 

\vspace{0.2cm}
\noindent 
Recall that $\gamma$ 
is oriented. It makes then sense to consider 
$L(\gamma)$ which is the set of punctures which are to 
the right of $\gamma$ but to the left of $[\gamma]$, and 
similarly $R(\gamma)$ which is the set of punctures which are to 
the left of $\gamma$ but to the right of $[\gamma]$. 
In order to define them properly let us consider 
the disks $D^+(\gamma)$ and $D^{-}(\gamma)$ bounded by 
$\gamma$ and arcs of the boundary of $P$. Note also 
by  $D^+([\gamma])$ and $D^{-}([\gamma])$ the respective disks 
in the case of $[\gamma]$. We assume that the positive 
disks lie on the left of the arc. Then 
$R(\gamma)$ is the set of punctures contained 
within  $D^+(\gamma)-D^-([\gamma])$ and 
$L(\gamma)$ the set of those from $D^-(\gamma)-D^+([\gamma])$. 
Notice that the punctures of $\gamma$ and 
$[\gamma]$  may be distinct. The punctures
of $L(\gamma)\cup R(\gamma)$ will be labeled $L$ and $R$ 
respectively.

\vspace{0.2cm}
\noindent
Remark that there exists a homeomorphism of  the punctured 
polygon $P^*$ which sends $\gamma$ onto $[\gamma]$. This 
implies that $L(\gamma)$ has the same cardinality as 
$R(\gamma)$.

\begin{center}
\includegraphics{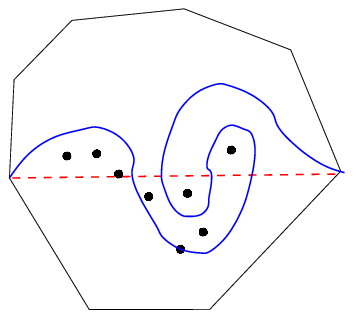}
\end{center}
Then the arc $\gamma$ can be seen as intersecting 
the ideal arcs in order (disregarding the position of 
the ideal arcs in the plane) some arcs 
might being crossed twice:

\begin{center}
\includegraphics{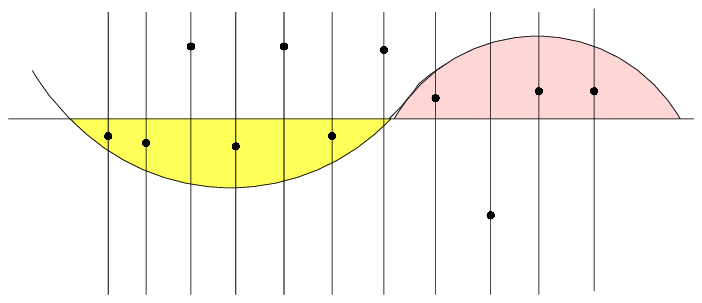}
\end{center}

\noindent 
The punctures on crossed arcs which do not belong to 
$L(\gamma)\cup R(\gamma)$ will be called {\em remote} punctures.

\vspace{0.2cm}
\noindent 
The main idea of the simplification procedure is that 
once we have two  consecutive  connected regions, one from 
$D^+(\gamma)-D^-([\gamma])$ and the next one 
from $D^-(\gamma)-D^+([\gamma])$, we are able to 
simplify the arc $\gamma$. In our simplification the 
line segment $[\gamma]$ will be made curvilinear, but we ignore 
this in our drawings.

\vspace{0.2cm}
\noindent Here is a sample: consider a pair of  
punctures having different labels which are closest to each 
other, thus the arc $\gamma$ does not meet but inert arcs 
in between. Assume that the first puncture is labeled $R$ and 
the second $L$. One might assume  that the inert punctures are far 
away in the plane. 
Consider the segment joining the two points and a 
small neighborhood of it, and further 
the braid twist which moves the two punctures clockwisely
by interchanging them. The effect of this move is as follows: 

\begin{center}
\includegraphics{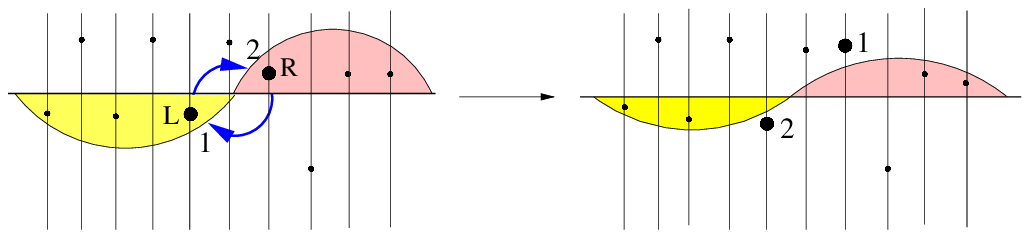}
\end{center}
In order to express the braid twist as a product of 
braid generators of $B_{\infty}$ one has to take care of the 
position of the inert punctures in the edges crossed by 
the segment relating the two punctures to be interchanged. 
If the left puncture is on the edge $e$ and the right one on the 
edge $f$ let $e_1,e_2,\ldots, e_p$ be the edges which are 
intersected by the arc $\gamma$. Since $e$ and $f$ are the closest 
possible then this sequence is made of distinct edges. 
We set then 
\[C_{ef}=\sigma_{ee_1}^{-1}\sigma_{e_1e_2}^{-1}\cdots 
\sigma_{e_{p-1}e_p}^{-1}\sigma_{e_{p}f}^{-1}
\sigma_{e_{p-1}e_p}\cdots \sigma_{e_{1}e_2}\sigma_{ee_1}\]
Notice that the edges $e_j$ and $e_{j+1}$ are adjacent edges 
of the triangulation $\tau_B^*$. The untangling factor 
$C_{ef}$ will reduce the number of  punctures in the two 
adjacent regions with different labels. 

\vspace{0.2cm}
\noindent
We continue to use inductive untangling terms for all punctures from $L$ 
sitting on the left of the puncture which switched from $L$ to $R$.
These have as effect the translation of that puncture on the leftmost 
side, where it becomes inert.  In particular the arc $\gamma$ 
becomes straight on the first quadrilateral that meets. 
Further we continue the same procedure for the next pair of 
punctures having different labels until one regions is empty. 
Remark that the procedure we describe is not canonical, other 
ways to switch the punctures might equally being used. 

\vspace{0.2cm}
\noindent
Let consider further the triangles $T_3,T_4,...$ crossed by 
the arc $\gamma$, which are not necessarily distinct.  Let also 
denote by $e=e_0,e_1,e_2,...$ the edges crossed by the 
arc $\gamma$ in the order. The arc $\gamma$ can cross once 
an edge $e_N$ at its midpoint puncture, in the case when the 
puncture associated to $\gamma$ is the same as the puncture on 
$e_N$. All the other intersections between $\gamma$ and the 
edges $e_j$ are different from the punctures, and thus it makes 
sense to say that $\gamma$ intersects $e_j$ to the 
left or to the right. Further $\gamma$ crosses $e_N$
both to the left and to the right. 

\vspace{0.2cm}
\noindent
Assume that $\gamma$ goes to the left and it is nonstraight. 
Let $e_j$ be the first edge so that $\gamma$ intersects 
$e_j$ to its right. According to our convention such a point 
exists always, but it might be 
the midpoint puncture of $\gamma$. We say that the punctures 
of $e$ and of $f=e_j$ are conjugate along  $\gamma$. 
Let $\sigma_{ef}$ be the braiding of the punctures 
of $e$ and $f$ along the arc $\gamma$. We wish to exchange these 
two punctures by means of a braiding. 

\vspace{0.2cm}
\noindent
Let $S_{ef}=C_{ee_1}C_ {e_1e_2}\cdots C_{e_{j-1}f}$. 
We call $S_{ef}$ the  first  untangling braid of the arc 
$\gamma$. By modifying the arc $\gamma$ by means of an untangling  
braid  we obtain an arc for which the cardinal of 
$L(\gamma)$  was decreased by one unit. Moreover, we have: 

\begin{lemma}
Let $\gamma'=S_{ef}\gamma$ be the image after 
the braid untangling. Then $\gamma'$ is straight with respect 
to $e$. 
\end{lemma} 
\begin{proof}
If $\gamma'$ were nonstraight then this would mean that 
the first puncture near which $\gamma'$ passed (namely 
the  central puncture of the first quadrilateral $Q$  which 
$\gamma'$ intersects)  would be in $L(\gamma)\cup R(\gamma)$. 
But we chose the untangling braid so that the central  
puncture in the first quadrilateral becomes remote. This implies  
the claim. 
\end{proof}

\vspace{0.2cm}
\noindent In particular, if  the arc $\gamma$ were not straight then 
we use first the braid untangling and further perform the flip like 
in Mosher's algorithm. We continue then with the same procedure 
until the arc is combed. When the arc is combed then we 
have to check whether it is straight. If not, then let us 
restart the straightening procedure once again and this time we use 
only  untangling braid  terms until we get stacked.  

\begin{lemma}
If $\gamma$ is a combed arc which admits no untangling braids
then $\gamma$ is straight.   
\end{lemma} 
\begin{proof}
This means that the sets $L(\gamma)$ and $R(\gamma)$ are empty, and 
thus we can use an isotopy keeping fixed the punctures which 
transforms $\gamma$ into its shadow line segment. 
\end{proof}

\vspace{0.2cm}
\noindent This means that eventually we transformed the arc 
$\gamma$ into an arc which belongs to $\tau_B^*$.

\subsection{Rectification of punctured triangulations}
In the process of combing an arc we  have to use flips which change the 
base (punctured) triangulation. In meantime the  untangling braid 
factors do not affect the reference triangulation but only the 
arc to be straightened. 
However, suitable sequences of flips could 
lead to edges which are tangled. We would like to keep the 
reference triangulations as simple as possible in the process of combing
in order to prevent them to have their edges too distorted. 
The way to do this is to rectify from time to time 
the triangulation. Notice that there is no analogous transformation 
in the case of the group $T$, since in that case  the 
flip of a geodesic triangulation is still a geodesic triangulation, 
as  two vertices define uniquely the segment joining them. 
In the punctured case we have to specify for any edge 
the corresponding puncture which belongs to it.

\vspace{0.2cm}
\noindent 
Let us give an example. It may happen that in the 
combing process we can obtain two different combinatorics of 
punctured triangulations of the pentagon, by choosing 
different locations for the interior punctures, as follows:

\begin{center}
\includegraphics{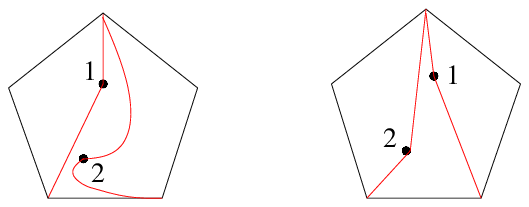}
\end{center}

\vspace{0.2cm}
\noindent We have then to choose one  triangulation which is allowed to 
appear in the combing process and discard the other one.  

\vspace{0.2cm}
\noindent
Specifically, 
for each triangulation $\tau$ (without punctures) 
we will choose a lift of it as a punctured 
triangulation $\tau^*$, which will be called {\em admissible}. 
If $\tau$ is supported in the polygon $P$ (thus it coincides with 
the Farey triangulation outside $P$) then we assume that $\tau_*$ has also 
the support in $P$. Moreover, we will ask that 
if $P\subset Q$ and 
$\tau_{Q}|_P =\tau_P$ then  $\tau^*_{Q}|_P =\tau^*_P$. 
Eventually, we can choose the set of admissible triangulations 
$\tau^*$ so that all of them can be obtained from each other 
by using flips or some braid.

\vspace{0.2cm}
\noindent It is convenient to choose the set of admissible 
triangulations of a pentagon as follows: 

\begin{center}
\includegraphics{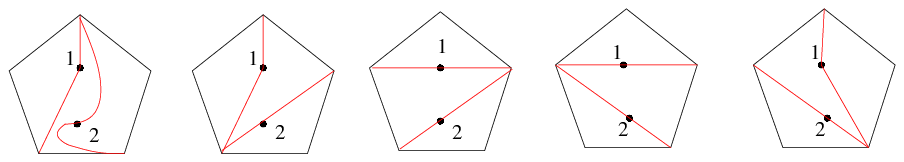}
\end{center}

\vspace{0.2cm}
\noindent Assume now that in the combing process at some step 
we modified the  basic punctured triangulation by a sequence at flips and get
$\delta^*$ which is not anymore admissible. Then there exists an unique 
admissible punctured triangulation $\tau^*$ which defines the same
subjacent triangulation as $\delta^*$ after dropping the punctures. 
There exists  then an unique braid  element in $B_{\infty}$ which transforms 
$\delta^*$ into $\tau^*$. We call it the {\em rectification (or correction) 
factor} and denote it ${\rm Cor}_{P}$ where $P$ is the 
finite support of the triangulation. 
For example, in the picture above the braid which interchanges 
counter-clockwisely the punctures $1$ and $2$ is the 
correction factor sending the left side pentagon into the right side 
pentagon.

\vspace{0.2cm}
\noindent 
We will assume that after any flip we will look at the newly obtained 
triangulation. If this is admissible then  we continue the straightening 
and combing. Otherwise, then we   
insert the correction factor which transform the intermediate 
triangulation  into an admissible punctured 
triangulation.

\subsection{The combing of $T^*$ is asynchronously bounded} 
\paragraph{Rewriting the raw combing into the two-generator
alphabet.} 
We proceed now in the same way as we did in the case of the group $T$.
Let $\zeta\in T^*$ be an element which is presented in the form of 
a couple of  punctured triangulations $(\tau_B^*,\tau_*)$. 
The arcs from $\tau_*$ coincide with those from $\tau_B^*$ outside 
some finite punctured polygon $P$. 
We have an order on the finite set of arcs 
$\gamma_1,...,\gamma_N$ of $\tau_*$ which belong to $P$,  
and each arc is oriented. 

\vspace{0.2cm}
\noindent
We start to comb and straight the arc $\gamma_1$. We record all 
transformations needed, namely the flips and the untangling  
braid terms, whose composition is $X_1$. 
Notice that these transformations do not affect the 
complementary of the support polygon $P$. One $\gamma_1$ is 
simplified and thus transformed into an arc $X_1(\gamma_1)$ of $\tau_B^*$ we 
look upon the image of  $\gamma_2$ under the transformation $X_1$. 

\vspace{0.2cm}
\noindent
Use next the combing and straightening procedure for $X_1(\gamma_2)$. 
Recall that arcs are disjoint and thus all elements used in 
combing and straightening (flips and braids) $X_1(\gamma_2)$ are supported 
in the complement of $X_1(\gamma_1)$, in particular they 
do not alter this arc. If the composition of these elements is 
$X_2$ then $X_2\circ X_1(\gamma_1)=X_1(\gamma_1), X_2\circ X_1(\gamma_2)$ 
belong to $\tau_B^*$. 

\vspace{0.2cm}
\noindent
Further continue the same procedure for all arcs in $P$. Notice 
that an arc $\gamma_k$ which belonged initially to $\tau_B^*$ 
will be left unchanged. This stops when all arcs from $P$ and thus 
all arcs from $\tau^*$ have been simplified to arcs in $\tau_B^*$. 

\vspace{0.2cm}
\noindent
The combing of $\zeta$ records the concatenation of sequences of flips,  
braid untanglings and correction factors 
used in $X_1,X_2,...,X_N$ in order to transform 
$\tau_*$ into $\tau_B^*$. This raw combing  of $T^*$  uses the 
infinite alphabet of all flips and all braid twists $\sigma_{ef}$, 
-  where $e$ and $f$ are adjacent edges -
from a suitable set of generators of $B_{\infty}$. Notice that 
the braid  terms $C_{pq}$ are canonically written as 
products of braid twists.  

\vspace{0.2cm}
\noindent
The second step consists of translation the raw combing into a combing 
based on the alphabet $\{\aps,\bps\}$, which is the generator system for 
$T^*$. One possibility is to use the normal forms determined 
in \cite{FK} for the braid elements as a words in $\aps,\bps$. 

\vspace{0.2cm}
\noindent
Another way to translate this is at follows.  Recall that 
the generator $\sigma_{ef}$ can be expressed in a very simple form, 
according to the results from \cite{FK2}. In fact,  it can be verified 
that the braid generator $\sigma_{[02]}$ associated to the edge joining the 
punctures numbered $0$ and $2$ of the support of $\aps$ can be written as  
\[ \sigma_{[02]}=(\bps\aps)^5\]
In other words, if $f$ is assumed to be the d.o.e. then 
$\sigma_{ef}$ is either  $(FR)^5$ or   $(RF)^5$, depending
on the relative position of $e$  with respect to $f$ 
(to its left side of to its right side, respectively).

\paragraph{Transfers for $T^*$.}
Further, in order to  get the word associated to  the composition 
of two braids one needs to insert  {\em transfers}, as we did when 
we translated the combing of $T$ in the two-generator alphabet 
(see section 2.4). Thus, when replacing in the 
product  $\sigma_{ef}\sigma_{fg}$  each  braid generator by the 
corresponding word in $F,R$ we have to insert in between the 
transfer $T_{gf}$ (also expressed as a word in $F$ and $R$). 
It suffices to find explicit formulas for the transfers $T_{fg}$. 
However, the same formulas that we used in section 3.4. for defining the 
transfer $T_{fg}$ as a word in $\alpha,\beta$ can be used now 
for the transfer $T_{fg}^*$, as a word in $\aps,\bps$. This is a 
consequence of the following splitting result:  

\begin{lemma}
The surjection $T^*\to T$ splits over $PSL(2,\Z)$ and thus we have a natural 
embedding $PSL(2,\Z)\subset  T^*$. 
\end{lemma}
\begin{proof}
Since $\aps^2$ is of order two and $\bps$ is of order 3 in $T^*$ and 
their  free amalgam $\langle\aps^2\rangle * \langle\bps\rangle$
(which is $PSL(2,\Z)$) embeds into $T$, then 
it will be also embedded in $T^*$, because $T^*$ surjects onto
$T$. Thus we also have a natural embedding of $PSL(2,\Z)$
into $T^*$ given by $\aps$ and $\bps$. 
\end{proof}
\noindent This procedure defines a combing for the group $T^*$ which 
uses the alphabet $\aps,\bps$.

\paragraph{$T^*$ is asynchronously combable.}
We need to prove that the combing of $T^*$  defined above is 
asynchronously bounded. Actually the proof given in section 2.5 
for $T$ can be adapted to $T^*$ with minor modifications. 
Since the form of the transfers is the same for $T^*$ and $T$, 
it suffices to look upon the way that the untangling braid and 
correction  factors 
intervene in the combings.

\vspace{0.2cm}
\noindent
Recall that we have to compare the path combings associated to nearby 
elements of $T^*$ and thus to elements that differ from each other 
by a factor $\aps$. Elements of $T^*$ are associated to triangulations. 
We have thus to make comparison between the combings simplifying 
the triangulations $\tau^*$ and $\alpha\tau^*$ in order to arrive at 
the same base triangulation $\tau_B^*$. Notice that this is the same 
to consider the pair of triangulations $(\tau_B^*,\tau_*)$ and 
$(\alpha^{-1}\tau_B^*,\tau^*)$. 
But this is the same as writing  the simplification procedure 
for a given triangulation $\tau^*$ with  
respect to two base triangulations $\tau_B^*$ and $F^{-1}\tau_B^*$. 
The advantage is that we can immediately see that the untangling braid
factors needed are  {\em almost} the same in the two cases 
and they depend mostly on  the topology of the arc in the 
complementary of the punctures.

\vspace{0.2cm}
\noindent 
We have then to understand the following situation: we have a sequences
of arcs $\gamma=\gamma_1, \gamma_2,...,\gamma_N$ that have to be straightened
and combed with respect first to $\tau_B^*$, and second to 
$F^{-1}\tau_B^*$. Then we have to compare the two combings and see whether 
they are at bounded asynchronous distance. 
We will analyze first the case of one arc $\gamma$ and observe that 
after straightening and combing it, then the situation for the next 
arc $\gamma_2$ is similar, namely we will have to comb/straighten $\gamma_2$ 
with respect to two base triangulations that differ from each other 
by a flip. Notice that  these new base triangulations have changed 
in meantime (several edges have been flipped in meantime and some 
braid corrections applied). 

\vspace{0.2cm}
\noindent 
Let thus concentrate on the case where we deal with one arc $\gamma$. 
Going back to the proof concerning $T$ we see that we shall understand 
what happens when the arc $\gamma$ enters the quadrilateral $Q_e$ 
(of $\tau_B^*$) having the edge $e$ as a diagonal, 
where $e$ is the d.o.e. Then $\gamma$ is issued from the vertex of 
a pentagon $P$ containing $Q_e$. Moreover, the other 
triangulation $F^{-1}\tau_B^*$ 
corresponds to choosing the d.o.e. $e^*$.

\vspace{0.2cm}
\noindent
The pentagon $P$ contains two punctures in its interior, say $p$ and $q$. 
The arc $\gamma$ could  first encircle the two punctures $p$ and 
$q$ a number of times and then exit along one side. Since the mapping class group of the twice punctured pentagon is $\Z$ it follows that in the interior 
of $P$ there is an unique configuration possible, namely: 
 
\begin{center}
\includegraphics{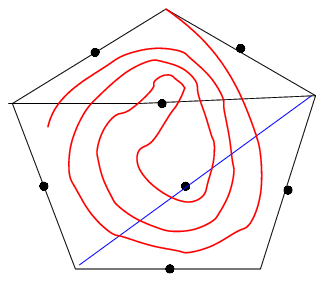}
\end{center}

\vspace{0.2cm}
\noindent Then one of the punctures $p$ and $q$ is conjugate to the 
first puncture that $\gamma$ meets and thus the first untangling 
braid factor is $\sigma_{pq}^k$ for some $k$. 
Thus we have to express this braid factor using the two possibilities for 
the d.o.e., namely $e$ and $e^*$. Using one base triangulation we 
have $\sigma_{pq}=(FR)^5$ while in the other one it reads 
$\sigma_{pq}=(RF)^5$. This means that the two words which describe the 
respective untangling braids are  
$(FR)^{5k}$ and respectively $F^{-1}(FR)^{5k}F$, which are  
at bounded asynchronous distance for any $k$.

\vspace{0.2cm}
\noindent 
After untangling the arc $\gamma$ we obtain an arc that encircles 
once  either one, or both or  else none of the two interior punctures  and 
then exits along one side. Notice that the arc $\gamma$ might 
return and run across $P$ again. 
There  are several cases to be taken into account (up to symmetry and 
the choice of the exit half-arc),    as follows:

\begin{center}
\includegraphics{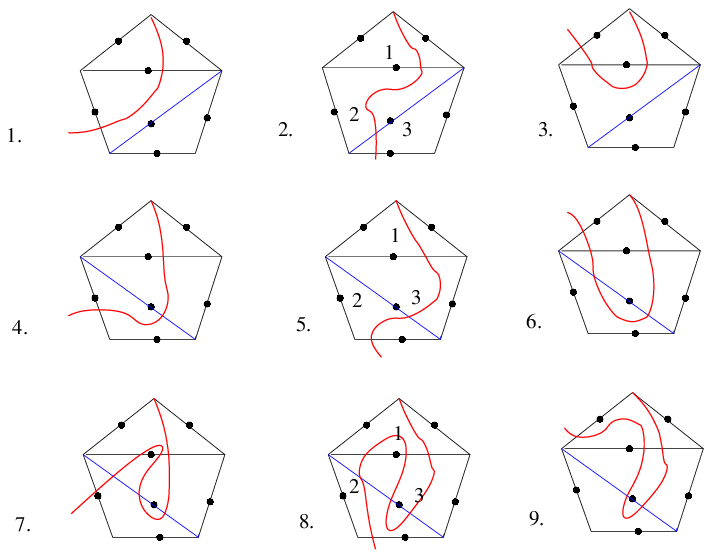}
\end{center}

\vspace{0.2cm}
\noindent 
Remark that the cases 7-9 are obtained by a braiding from 
1-3 and this braiding is the untangling factor in both situations, 
so that it suffices to check 1-6.  

\vspace{0.2cm}
\noindent 
If the arc  $\gamma$ does not intersect
neither $e$ nor $e^*$, or both of them, 
then the untangling factors and corrections factors
should be the same in both configurations.

\vspace{0.2cm}
\noindent 
Let us analyze these cases one by one. We will describe the 
simplification procedure with respect to $\tau_B^*$ and then 
$F^{-1}\tau_B^*$ at each step. 

\begin{enumerate}
\item 
\begin{enumerate} 
\item Assume that the arc $\gamma$ exits the pentagon $P$ after 
passing near the puncture 1 then encounters the punctures 2,3,...,m, 
and finally the puncture $r$ which is conjugate to $p$. 
We suppose for the moment that $\gamma$ does not return to the 
pentagon $P$. The straightening algorithm then for $\gamma$ 
with respect to $\tau_B^*$ runs at follows.
The first edge met by  $\gamma$ is $e_p$ and this cannot be flipped. 
We compute the untangling braid $C_{pr}$ which is the product of 
consecutive braid generators 
$\sigma_{rm}\sigma_{mm-1}\cdots \sigma_{21}\sigma_{1p}$. 
Untangle $\gamma$ in order to be able to flip $e_p$ and do the flip. 
The next edge is $e_1$  the edge containing $1$ and we continue.

\begin{center}
\includegraphics{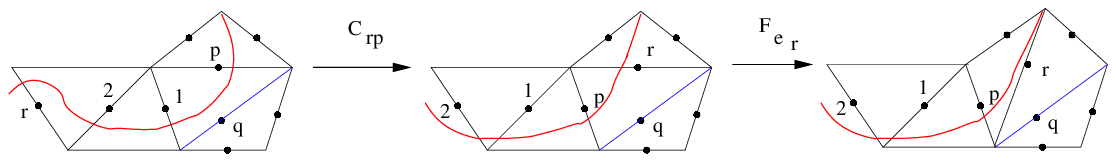}
\end{center}

\vspace{0.2cm}
\noindent 
\item 
If we deal with the d.o.e. $e^*$ then the straightening process is different. 
The first edge to be combed is $e_p$ but the second one is $e_q$ instead of 
$e_1$ which comes in third. Moreover, the first puncture conjugate 
to $p$ is not anymore $r$ but $q$. We have then to untangle the arc 
$\gamma$ by using $C_{pq}$ and then we can flip the first edge 
met by $\gamma$. The next edge to be flipped is $e_1$. Moreover the 
puncture $p$ (which is now located where $q$ lived before) is conjugate 
to the same $r$ as in the previous situation. 
We use then the untangling factor $C_{pr}$ and then flip 
the edge $e_r$. We obtained two triangulations on the pentagon $P$, that 
obtained in the figure above and the present one. 
We assumed that the former is admissible  and so  
our present triangulation is not admissible and it has to be 
corrected by a factor ${\rm Cor}_{P}$. The correction factor is actually 
the braid generator $\sigma_{qr}^{-1}$ in this case. 

\begin{center}
\includegraphics{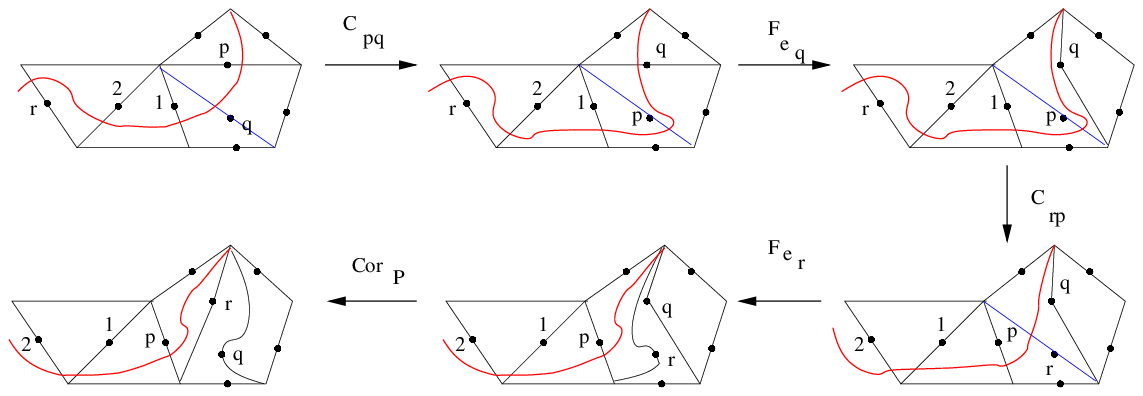}
\end{center}

\end{enumerate}

\vspace{0.2cm}
\noindent 
The output of the two simplification procedures is almost the same 
in the two situations, namely the two final diagrams differ by only one flip 
on the edge $F_{e_q}$ located within the pentagon $P$. 
In particular the procedure will be from now on the same in both situations. 
Since we supposed that the arc does not renter $P$ then we get 
asynchronous boundedness on the second part.

\vspace{0.2cm}
\noindent The only trouble we had in the first part is that 
the puncture $r$ might be far away, and thus we have to understand 
the distance between the combing when using the untangling procedure. 
However,  the untangling factors are the same except for their 
last terms, and they are products of braid generators. 

\vspace{0.2cm}
\noindent 
In order to apply a braid generator, say $\sigma_{rm}$ we have to use 
the transfer of the d.o.e. $e$ to the position $e_m$, then insert the 
classical braid generator $(RF)^5$ (or $(FR)^5$) and further 
come back using the inverse transfer $T_{e_me}$.
But the same procedure was used for dealing with flips instead 
of braids when we described the choice of transfer making 
far away flips be at  asynchronous bounded distance in the two 
situations.

\vspace{0.2cm}
\noindent Recall that each untangling factor $C_{pq}$ is a product of 
several copies  of the standard braid generator $\sigma_{ef}=(RF)^5$
with the transfers $T_{fg}$. Moreover, the transfers used for the 
punctured case coincide with those for the group $T$. Thus, the 
two combings of nearby elements are obtained from the combings in $T$ 
by inserting finitely many elements $(RF)^5$. In particular, since the 
combing of $T$ is asynchronously bounded and we insert  
only elements of bounded length  we obtain path at asynchronously bounded 
length.

\vspace{0.2cm}
\noindent The second problem that we might encounter is that 
the arc $\gamma$ might return within $P$ before reaching the 
first puncture conjugate to $p$. According to our disk lemma then 
$\gamma$ will renter the edge $e_1$ on the other side of the 
puncture $e_1$. However, each new entrance in $P$ will 
contribute to the first untangling factor by a braid element which 
interchanges two punctures inside the pentagon. The words 
which describe the braid elements  
with respect to the two triangulations are conjugate to each other 
by a factor $F$. Therefore there is no problem in keeping them 
at asynchronous bounded distance.  

\vspace{0.2cm}
\noindent For the remaining cases we will just picture the 
simplification steps within the pentagon $P$, until the moment from 
where the reductions are the same in both situations. 

\item 
\begin{enumerate}
\item Here is the first simplification for 2: 

\begin{center}
\includegraphics{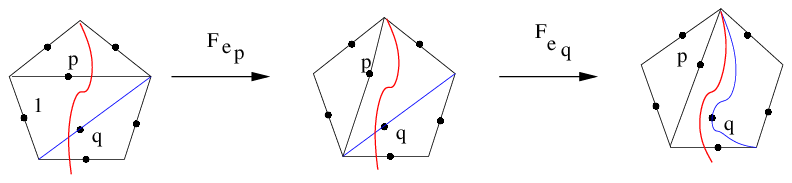}
\end{center}

\item The simplification with respect to $F^{-1}\tau_B^*$:

\begin{center}
\includegraphics{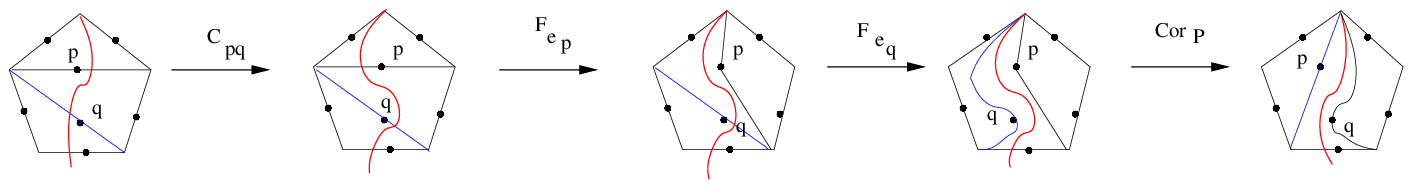}
\end{center}

\end{enumerate}

\item In the case 3 the arc $\gamma$ does not intersect neither 
$e$ nor $e^*$. The simplifications will be therefore 
identical in the two situations. 

\item The case 4. We assume that $q$ is conjugate to some puncture $s$. 

\begin{enumerate}
\item Here is the first simplification for 4: 

\begin{center}
\includegraphics{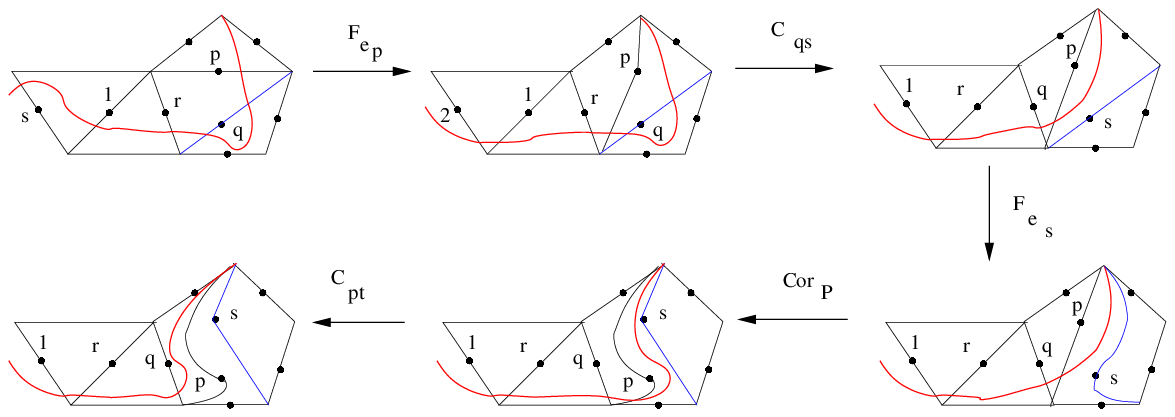}
\end{center}

\noindent where $t$ is the puncture conjugate to $p$ at the fifth 
stage. 

\item The simplification with respect to $F^{-1}\tau_B^*$:

\begin{center}
\includegraphics{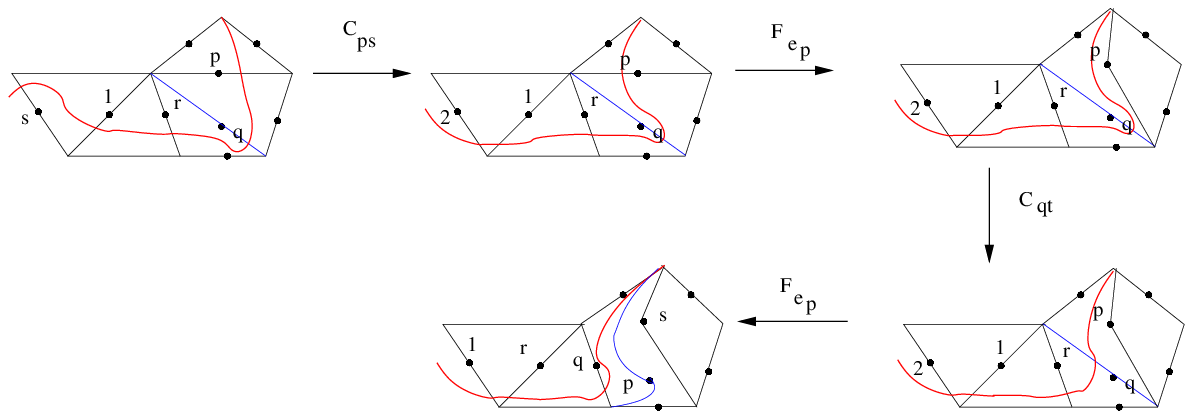}
\end{center}

\noindent Remark that at the third step the puncture $q$ is conjugate to the 
same puncture $t$ which appeared above. 
\end{enumerate}

\item The case 5. 

\begin{enumerate}
\item Here is the first simplification for 5: 

\begin{center}
\includegraphics{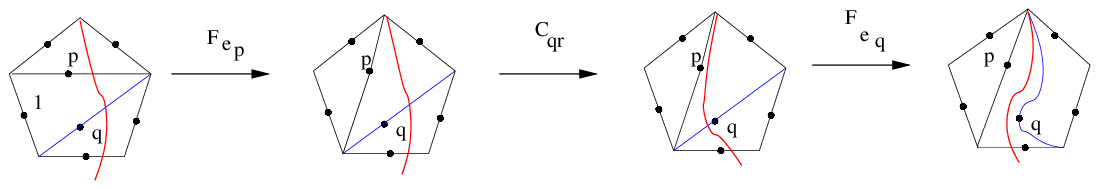}
\end{center}

\item The simplification with respect to $F^{-1}\tau_B^*$:

\begin{center}
\includegraphics{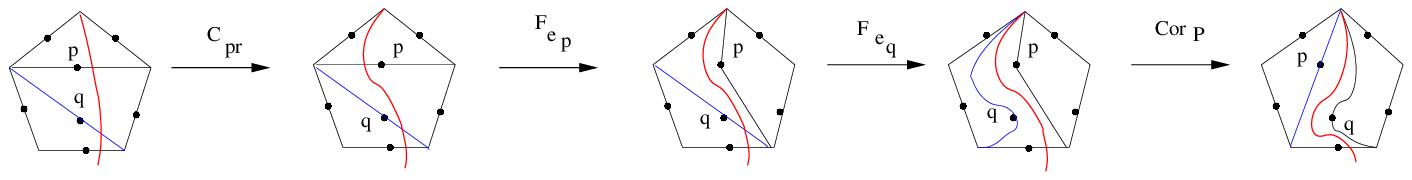}
\end{center}

\noindent Remark that  $p$ is conjugate to the 
same puncture $r$ that appeared above. 
\end{enumerate}

\item Last  simplification for 6: 
\begin{enumerate} 
\item The first situation is below: 

\begin{center}
\includegraphics{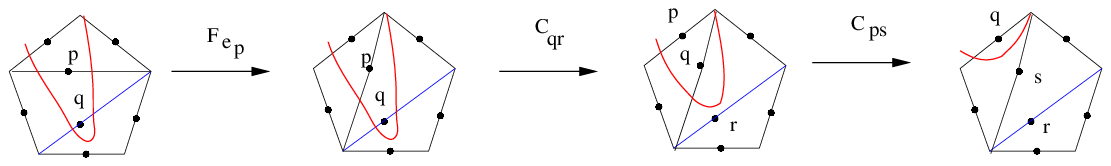}
\end{center}

\item The simplification with respect to $F^{-1}\tau_B^*$:

\begin{center}
\includegraphics{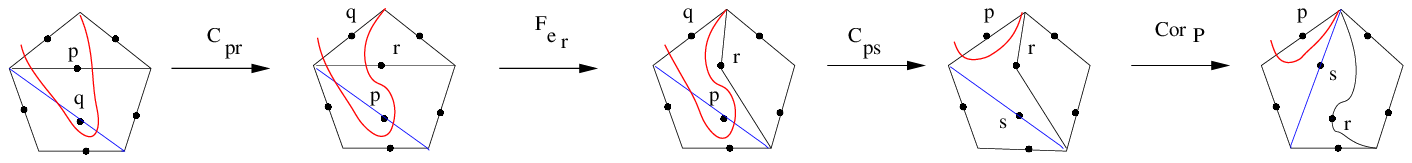}
\end{center}

\end{enumerate}

\end{enumerate}

\vspace{0.2cm}
\noindent This analysis shows that the combing that we defined 
is asynchronously bounded, as claimed.

\subsection{The departure function}

The combing that we defined has not a departure function. 
In fact, let us assume that the Mosher combing consists of two 
flips which are located on nearby edges $f_1,f_2$, which are faraway 
from the d.o.e. $e$. Our procedure amounted to translate the  
sequence $F_{f_1}, F_{f_{2}}$ into 
$T_{e, f_1}, F, T_{f_1^*, e}, T_{e,f_2}, F, T_{f_2^*,e}$. 
However, if $f_1$ and $f_2$ are nearby edges then the transfers 
$T_{e,f_1^*}$ and $T_{e,f_2}$  are paths having  in common a  
large part of their initial segments. Furthermore, the 
modified transfers will be also kept at finite distance one from the 
other along this initial segment.  This will contradict the 
departure function condition since one finds a long path having 
the shape of a back and forth path in the Cayley graph.    

\vspace{0.2cm}
\noindent 
This is however the only accident that might occur. 
In fact let us assume that whenever we have two  consecutive 
transfers $T_{f_1,e} T_{e,f_2}$ in the translation of 
Mosher's combing then we simplify it. This means that 
we consider the geodesics joining $e$ to $f_1$ and respectively $f_2$ 
and drop their common part. We obtain then a word which is actually 
the representative of  $T_{f_1,f_2}$, and whose   
terms of the form $aa^{-1}$ have been cancelled.  
The combing so obtained is called the reduced combing of $T^*$.

\begin{proposition}
The reduced combing is still asynchronously combable and moreover 
it has a departure function. 
\end{proposition}
\begin{proof}
The reduced combing just cancelled the terms  of the form $aa^{-1}$ and so the asynchronous combability is preserved. 

\vspace{0.2cm}
\noindent 
We claim first that if $\tau_1$ and $\tau_2$ agree 
outside the ball $B(r)$ of radius $r$ (of the dual tree) 
and their d.o.e. lay within 
$B(r)$ then any time that $\tau_1,\ldots, \mu,\ldots,\tau_2$ 
appears as a chain in the combing of some element of $T$ then 
any intermediary triangulation $\mu$ should agree  with 
$\tau_i$ outside the ball of radius $r$. 
Since there are finitely many 
triangulations of the polygon with $3\cdot 2^{r}$ sides (corresponding to 
the ball $B(r)$ by duality) we find that 
the length of such a chain is bounded by a function on $r$. 
In particular if the length of such a chain is bigger than that function 
then $\tau_i$ should be distinct outside the ball 
of radius $r$. However if we want to modify the standard triangulation 
by means of an element of $T$ written as a word $w$ in $\alpha,\beta$ 
then the newly obtained triangulation coincide with the former 
outside the ball of radius $|w|+1$, where $|w|$ denotes the 
length of the word $w$. In particular the element 
$(\tau_1,\tau_2)\in T$ should have length at least $r-1$ and 
so the elements that $(\tau,\tau_1)$ and $(\tau,\tau_2)$ are at 
distance $r-1$ far apart in the Cayley graph. 
This will prove therefore that the  reduced combing of $T$ has a 
departure function. 

\vspace{0.2cm}
\noindent 
The proof of our claim is at follows. Each chain in a combing is made 
of transfers and flips. If there exists a transfer going outside 
the ball of radius $B(r)$ then either there exists a flip at the end 
of the transfer or else there exists a transfer going back. 
However any time that we used a flip, this is done for combing, and thus 
the edge remains flipped until the combing terminates. 
Thus $\tau_2$ should contain the edge flipped outside $B(r)$ which is 
a contradiction. The other possibility was to have a transfer that goes out 
of $B(r)$ and then enters again $B(r)$, but this is impossible since 
the ball  in the tree are convex with respect to the (modified) 
transfers. This means that two edges 
$e,f$ lying inside $B(r)$ are joined by a  modified transfer 
which does not affect the complementary of $B(r)$. 
The last possibility is that we have a composition of two transfers 
that go outside the ball and then come back. 
This is impossible since then we should have 
a common part in the two transfers, but 
our reduction procedure cancelled such terms.

\vspace{0.2cm}
\noindent The claim holds true also for punctured triangulations, when 
dealing with the group $T^*$. However it does not imply directly 
the existence of the departure function since the number of 
punctured triangulations of a given polygon is infinite. 
In meantime our hypothesis should be stronger. 
The punctured triangulations are not only agreeing outside 
$B(r)$ but also inside the ball $B(r)$ they do not differ too much from 
each other, in order to be close in the Cayley graph. 
This means that  one obtains $\tau_2$ from $\tau_1$ by means 
of a sequence of abstract combings and a braid action, where 
the braid is viewed as an element of $T^*$ of bounded length 
(say $N$) as a word in $\alpha^*,\beta^*$.   
Since the number of abstract triangulations of the polygon with  $3\cdot 2^r$
sides is finite and the length of the braid is less than $N$ 
we find that the length of a chain joining $\tau_1$ to $\tau_2$ 
in a combing process is bounded in terms of $r, N$.

\end{proof}

{\small     
\bibliographystyle{plain}      
      
}      
\end{document}